\magnification=\magstephalf
\input eplain
\input BMmacs

\phantom{x}\vskip .5cm
\centerline{{\titlefont ON THE CHOW MOTIVE OF AN ABELIAN SCHEME}}
\medskip

\centerline{{\titlefont WITH NON-TRIVIAL ENDOMORPHISMS}}
\vskip .5cm

\centerline{{\byfont by}\quad {\namefont Ben Moonen}}
\vskip 1cm

{\eightpoint 
\noindent
{\bf Abstract.} Let $X$ be an abelian scheme over a base variety~$S$ with endomorphism algebra~$D$. We prove that the relative Chow motive $R(X/S)$ has a canonical decomposition as a direct sum of motives~$R^{(\xi)}$ where $\xi$ runs over an explicitly determined finite set of irreducible representations of the group $D^{\opp,*}$, such that $R^{(\xi)}$, seen as a functor from Chow motives to $D^{\opp,*}$-representations, is $\xi$-isotypic. Our decomposition refines the motivic decomposition of Deninger and Murre, as well as Beauville's decomposition of the Chow group. 

The second main result is that we construct a canonical generalized motivic Lefschetz decomposition. Inspired by work of Looijenga and Lunts in~[\ref{LooijLunts}], the role of ${\frak sl}_2$ in the classical theory is here replaced by a larger Lie algebra, generated by all Lefschetz and Lambda operators associated to non-degenerate line bundles. We construct an action of this larger Lie algebra on~$R(X/S)$ and deduce from this a generalized Lefschetz decomposition. We also give a precise structure result for the Lefschetz components that arise.

As an application of our techniques, we provide a positive answer to a question of Claire Voisin concerning the analogue for abelian varieties of a conjecture of Beauville.
\smallskip

\par}
\vskip 8mm

\tocheader
\readtocfile
\par

\intro

\ssection
As an application of Fourier theory, Beauville proved in~[\ref{BeauvChow}] that the Chow ring (with $\mQ$-coefficients) of a $g$-dimensional abelian variety~$X$ has a bigrading $\CH(X) = \oplus_{j,s}\, \CH^j_{(s)}(X)$, where the upper grading is given by the codimension of cycles and $[m]_X^*$ acts on $\CH^j_{(s)}(X)$ as multiplication by $m^{2j-s}$. As shown by Deninger and Murre in~[\ref{DenMur}], this decomposition in fact comes from a natural decomposition $R(X) = \oplus_{i=0}^{2g}\, R^i(X)$ of the Chow motive of~$X$; we have $\CH^j_{(s)}(X) = \CH^j\bigl(R^{2j-s}(X)\bigr)$. The results of Deninger and Murre are valid, more generally, for abelian schemes $X \to S$ over a smooth quasi-projective base variety over a field. 

One way to state Beauville's result is by saying that $\mQ^*$ acts on the Chow ring (letting $m/n \in \mQ^*$ act as $[m]_X^* \circ [n]_X^{*,-1}$), and that the only characters that occur in this representation are the characters $q \mapsto q^i$ for $i \in \{0,1,\ldots,2g\}$. The first main purpose of this paper is to explain how this can be refined in the presence of non-trivial endomorphisms.

To describe the result, consider an abelian scheme $X \to S$ of relative dimension~$g$ whose endomorphism algebra $D = \End^0(X/S)$ is  a simple algebra. (This is the essential case, to which the general case is reduced; see~\ref{GeneralCase}.) The group $D^{\opp,*}$ acts on~$\CH(X)$ and on the motives $R^i(X/S)$, which are objects of the category~$\Mot^0(S)$ of relative Chow motives over~$S$. This induces the structure of a $D^{\opp,*}$-representation on $\Hom_{\Mot^0(S)}\bigl(M,R(X/S)\bigr)$, for any relative Chow motive~$M$.

Let $G$ be $D^{\opp,*}$, viewed as a reductive group over~$\mQ$. The irreducible representations of~$G$ over~$\mQ$ are indexed by the $\Gamma$-orbits in a space~$\chargp^+$ of highest weight vectors, where $\Gamma$ is the absolute Galois group of the center of~$D$.  There is a natural ``weight function'' $\wt{\ } \colon \chargp^+/\Gamma \to \mZ$. We define an explicit finite subset $\chargp^\adm/\Gamma \subset \chargp^+/\Gamma$ of ``admissible'' elements, and we describe an involution $\xi \mapsto \xi^\wdual$ of this set, with $\wt{\xi^\wdual} = 2g-\wt{\xi}$; see the body of the text for the details. Our first main result, Theorem~\ref{MotDecRefined}, can then be stated as follows.

\ssectionnonr
{\it Theorem 1. --- There is a unique motivic decomposition
$$
R(X/S) = \bigoplus_{\xi \in \chargp^\adm/\Gamma}\, R^{(\xi)}(X/S) \leqno{(1)}
$$
that is stable under the action of~$D^{\opp,*}$ and has the property that for any motive~$M$ the $D^{\opp,*}$-representation $\Hom_{\Mot^0(S)}\bigl(M,R^{(\xi)}(X/S)\bigr)$ is $\xi$-isotypic. For $\xi \in \chargp^\adm/\Gamma$ we have $0 \leq \wt{\xi} \leq 2g$ and $R^i(X/S)$ is the direct sum of the motives~$R^{(\xi)}(X/S)$ with $\wt{\xi} = i$. Further we have a motivic Poincar\'e duality isomorphism $R^{(\xi)}(X/S)^\vee \isomarrow R^{(\xi^\wdual)}(X/S)\bigl(g\bigr)$ and, for $\xi \in \chargp^\adm/\Gamma$ with $\wt{\xi} = i$, a motivic Fourier duality $\Four\colon R^{(\xi)}(X/S) \isomarrow R^{(\xi^\wdual)}(X^\avdual/S)\bigl(g-i\bigr)$.\/}
\medskip

In particular, the Chow groups $\CH^j\bigl(R^{(\xi)}(X/S)\bigr) = \Hom_{\Mot^0(S)}\bigl(\unitmot(-j),R^{(\xi)}(X/S)\bigr)$ are $\xi$-isotypic as representations of~$D^{\opp,*}$.

\ssection
The second goal of the paper is to construct a canonical motivic Lefschetz decomposition. As K\"unnemann has shown in~[\ref{KunLef}], the choice of a polarization of~$X$ gives a motivic Lefschetz decomposition of $R(X/S)$ in the category $\Mot(S)$ of relative Chow motives with respect to ungraded correspondences. This decomposition depends on the choice of the polarization; it is not, in general, compatible with the decomposition of Theorem~1.

In [\ref{LooijLunts}] Looijenga and Lunts introduced the idea that one may replace the Lie algebra $\sl_2$ that gives the Lefschetz decomposition, by a larger ``N\'eron-Severi Lie algebra'' $\gg_\NS$ that is generated by the Lefschetz and Lambda operators associated to all possible choices of a polarization; they studied these operators on the cohomology of complex varieties and obtained some important results on the structure of the resulting Lie algebra.

Our canonical (generalized) Lefschetz decomposition of $R(X/S)$ is obtained from the action of a Lie algebra $\gg = \sp(X \times_S X^\avdual)$ on $R(X/S)$; see Theorem~\ref{RefLefThm}. This Lie algebra~$\gg$ is essentially the Lie algebra~$\gg_\NS$ of Looijenga and Lunts. (In some cases it is slightly bigger.) It has a grading $\gg = \gg_{-2} \oplus \gg_0 \oplus \gg_2$. The degree zero part~$\gg_0$ is given by the endomorphisms of~$X$. The elements of~$\gg_2$ (resp.\ $\gg_{-2}$) are the symmetric homomorphisms $X \to X^\avdual$ (resp.\ $X^\avdual\to X$), acting on $R(X/S)$ through Lefschetz operators (resp.\ Lambda operators) as defined by K\"unnemann. Our main new contribution is that we establish a commutation relation between arbitrary Lefschetz and Lambda operators. (In~[\ref{KunLef}] this is only done in the case where $L$ and $\Lambda$ come from the same polarization of~$X$.) We should note that the action of~$\gg$ on the motive is also closely related to earlier work of Mukai, Orlov and, especially, Polishchuk; see the introduction of Section~\ref{LieAction} for further discussion.   

As an application of Theorem~\ref{RefLefThm}, we obtain a generalized motivic Lefschetz decomposition. The result can be summarized as follows.

\ssectionnonr
{\it Theorem 2. --- There is a unique motivic decomposition
$$
R(X/S) = \bigoplus_{\psi \in \Irrep(\gg)}\, R_\psi(X/S) \leqno{(2)}
$$
in the category $\Mot(S)$ that is stable under the action of~$\gg$ and has the property that for any motive~$M$ the $\gg$-representation $\Hom_{\Mot(S)}\bigl(M,R_\psi(X/S)\bigr)$ is $\psi$-isotypic. This decomposition is stable under the action of~$D^{\opp,*}$ and there is a unique decomposition
$$
R(X/S) = \bigoplus_{\psi \in \Irrep(\gg)\atop \xi \in \chargp^\adm/\Gamma} \, R_\psi^{(\xi)}(X/S)
$$
in $\Mot^0(S)$ that is a common refinement of the decompositions\/~{\rm (1)} and\/~{\rm (2)}.}

\ssection
We refer to the motives $R_\psi(X/S)$ as the {\it Lefschetz components\/} of $R(X/S)$. There are only finitely many non-zero Lefschetz components; see Remark~\ref{HighWeightsRem} for how to determine which $\psi \in \Irrep(\gg)$ occur. 

For a Lefschetz component $R_\psi = R_\psi(X/S)$ there is an integer~$m$ such that $R_\psi$ lives in degree $g-m,g-m+2,\ldots,g+m$. The summand $R_\psi^\prim = R_\psi^{g-m}$ in the lowest degree is called the {\it primitive part\/} of~$R_\psi$. In the classical case, when $D=\mQ$, the Lefschetz component $R_\psi$ is isomorphic to $\oplus_{j=0}^m\, R_\psi^\prim(-j)$. In general, however, $R_\psi$ is no longer a sum of copies of Tate twists of~$R_\psi^\prim$ and something new happens. The primitive part should, in fact, not be thought of as a ``building block'' for the whole Lefschetz component. This role is taken by a smaller piece~$P_\psi(X/S)$ that we call the {\it core\/} of the Lefschetz motive.

\ssectionnonr
{\it Theorem 3. --- Let $\psi \colon \gg \to \gl(V)$ be an irreducible representation over~$\mQ$ such that the corresponding Lefschetz component $R_\psi(X/S)$ is non-zero, and consider the division algebra $B = \End_\gg(V)$. Then there exists a motive $P_\psi(X/S)$ with right $B$-action such that $R_\psi(X/S) \cong P_\psi(X/S) \otimes_B V$.}
\medskip

See Theorem~\ref{LefCompStr} in the text for a more precise version of this result, containing a more intrinsic description of the core. Let us note that the description of the Lefschetz component in the form ``$\hbox{core} \otimes_B \hbox{representation}$'' seems new even if we pass to realizations. The primitive part $R_\psi^\prim$ of a Lefschetz component is isomorphic to a sum of copies of the core~$P_\psi$.

\ssection
In the final section of the paper, we apply our techniques to answer a question of Claire Voisin. If $X$ is an abelian variety over a field we have Beauville's decomposition of $\CH(X)$, and it seems plausible that the associated descending filtration is a Bloch-Beilinson filtration. This speculation leads to some concrete, but very difficult, questions. For instance, it is expected that the cycle map should be injective on the zeroth layer $\CH_{(0)} = \oplus_j\, \CH^j_{(0)}(X)$, but this seems to be beyond reach of the presently available techniques. Voisin's question is whether the cycle class map is injective on the classes that are polynomial expressions in the symmetric divisor classes. This is the analogue for abelian varieties of Beauville's conjecture in~[\ref{BeauvSplitting}] about the so-called weak splitting property. We obtain a positive answer to Voisin's question.

\ssectionnonr
{\it Theorem 4. --- Let $X$ be an abelian variety over a field~$F$. Let $y \in \CH(X)$ be an element that can be written as $y = P(\ell_1,\ldots,\ell_r)$ for some polynomial $P \in \mQ[t_1,\ldots,t_r]$ and classes $\ell_i \in \CH^1_{(0)}(X)$. If $y$ is numerically trivial, $y=0$.\/}
\medskip

We have a similar result about the first ``layer'' $\CH_{(1)}(X)$; in this case we prove that the Abel-Jacobi map is injective on the set of classes $y \in \CH_{(1)}(X)$ that can be written as a polynomial in divisor classes. See Theorem~\ref{AJInjectThm}.

To conclude this introduction let us note that G.~Ancona, in his Paris thesis~[\ref{Ancona}], has independently obtained related results for universal abelian schemes over PEL Shimura varieties over~$\mC$.
\medskip

\ssectionnonr
{\it Acknowledgements.\/} --- I thank Christopher Deninger, H\'el\`ene Esnault, Klaus K\"unnemann and Frans Oort for their comments on an earlier version of the paper. I thank Claire Voisin for suggesting the problem discussed in Section~\ref{Applic1}.

\ssectionnonr
{\it Conventions.\/} --- Throughout, Chow groups are taken with $\mQ$-coefficients. Unless specified otherwise, all actions of groups or algebras are left actions. If $x$ is an element in a $\mQ$-algebra, we use the divided power notation $x^{[n]}$ for $x^n/n!$.

\section{Some inputs from representation theory}{RepTheory}

\ssection
\ssectlabel{BKkNotation}
In this section we consider a simple algebra~$B$ of finite dimension over a field~$k$ of characteristic~$0$. Let $K$ be the center of~$B$, let $[K:k] = n$ and $d = \dim_K(B)^{1/2}$. 

Let $\kbar$ be an algebraic closure of~$k$ and let $\Sigma(K)$ denote the set of $k$-algebra homomorphisms $K \to \kbar$. Let $\Ktilde$ denote the normal closure of~$K$ inside~$\kbar$, and write $\Gamma = \Gal(\Ktilde/k)$. The natural action of $\Gal(\kbar/k)$ on~$\Sigma(K)$ factors through an action of~$\Gamma$.

\ssection
\ssectlabel{HRootDatum}
Let $H$ be the reductive group over~$K$ with $H(R) = (B \otimes_K R)^*$ for any commutative $K$-algebra~$R$. Let $\bigl(\chargp(H), \Phi, \chargp^\vee(H), \Phi^\vee, \Delta\bigr)$ be the based root datum of~$H$. We need to recall the definition of~$\chargp(H)$; see for instance [\ref{WedhOrdin}], Section~1.2, for further details. Consider pairs $(T,Q)$ consisting of a maximal torus $T \subset H_\Kbar$ and a Borel subgroup $Q \subset H_\Kbar$ containing~$T$. Given such a pair, let $\chargp_{(T,Q)}$ denote the character group of~$T$. If $(T^\prime,Q^\prime)$ is another pair, there exists an element $h \in H(\Kbar)$ such that $hTh^{-1} = T^\prime$ and $hQh^{-1} = Q^\prime$. The induced isomorphism $\chargp_{(T^\prime,Q^\prime)} \isomarrow \chargp_{(T,Q)}$ is independent of the choice of~$h$ and $\chargp(H)$ is defined as the projective limit of the groups $\chargp_{(T,Q)}$. For any pair $(T,Q)$ the natural map $\chargp(H) \to \chargp_{(T,Q)}$ is an isomorphism.

There is a natural choice for an ordered $\mZ$-basis $\{e_1,\ldots,e_d\}$ of~$\chargp(H)$, obtained in the following way. Choose an isomorphism of $\Kbar$-algebras $a\colon B \otimes_K \Kbar \isomarrow M_d(\Kbar)$; this induces an isomorphism $\xi \colon H_\Kbar \isomarrow \GL_{d,\Kbar}$. Let $T \subset Q \subset H_\Kbar$ be the maximal torus and Borel subgroup such that $\xi(T)$ is the diagonal torus and $\xi(Q)$ is the upper triangular Borel. Let $\epsilon_j^\prime \colon \xi(T) \to \mG_{m,\Kbar}$ be the character that sends a diagonal matrix with entries $(c_1,\ldots,c_d)$ to~$c_j$, and define $\epsilon_j \in \chargp_{(T,Q)}$ by $\epsilon_j = \epsilon_j^\prime \circ \xi$. Then $\{\epsilon_1,\ldots,\epsilon_d\}$ is an ordered $\mZ$-basis of~$\chargp_{(T,Q)}$. Now define $\{e_1,\ldots,e_d\}$ to be the ordered $\mZ$-basis of~$\chargp(H)$ such that $e_j \mapsto \epsilon_j$ under the isomorphism $\chargp(H) \isomarrow \chargp_{(T,Q)}$. It follows from the Skolem-Noether theorem and the definition of~$\chargp(H)$ that the ordered basis thus obtained does not depend on the choice of the isomorphism~$a$. Further it is clear from the construction that the roots are the vectors $e_i - e_j$ for $i \neq j$, and that the basis of positive roots is given by $\Delta = \bigl\{e_i - e_{i+1} \bigm| i=1,\ldots,d-1\bigr\}$.

\ssection
\ssectlabel{HIrreps}
The group~$H$ is an inner form of $\GL_d$; hence the Galois group $\Gal(\ol{K}/K)$ acts trivially on the root datum of~$H$. By [\ref{Tits}], Thm.~7.2, we have a bijective correspondence between the set of irreducible finite-dimensional representations of~$H$ over~$K$ and the set $\chargp(H)^+$ of dominant weights. 

With respect to the ordered basis $\{e_1,\ldots,e_d\}$ as in~\ref{HRootDatum}, the dominant weights are the vectors $\lambda_1 e_1 + \cdots + \lambda_d e_d$ for $\lambda = (\lambda_1,\ldots,\lambda_d) \in \mZ^d$ with $\lambda_1 \geq \lambda_2 \geq \cdots \geq \lambda_d$. This gives an identification of~$\chargp(H)^+$ with the set
$$
\Lambda^+ = \bigl\{\lambda = (\lambda_1,\ldots,\lambda_d) \in \mZ^d \bigm| \lambda_1 \geq \lambda_2 \geq \cdots \geq \lambda_d\bigr\}\, . \eqlabel{Lambda+Def}
$$
For $\lambda \in \Lambda^+ = \chargp(H)^+$, let $\psi_\lambda$ be the corresponding irreducible representation of~$H$ over~$K$. 

If $\phi_\lambda$ is the irreducible representation of $\GL_d$ with highest weight given by~$\lambda$, the representation~$\psi_\lambda$ is a $K$-form of the representation $\phi_\lambda^{\oplus d(\lambda)}$ for some integer $d(\lambda)$ that divides~$d$. For later use, let us also recall that if $\lambda_d \geq 0$, the representation~$\phi_\lambda$ is the one obtained from the standard representation of~$\GL_d$ applying the Schur functor~$\mS_\lambda$. In the general case, without the assumption that $\lambda_d \geq 0$, we take an integer~$m$ with $\lambda_d + m \geq 0$; then $\phi_\lambda = \phi_{(\lambda_1+m,\ldots,\lambda_d+m)} \otimes \det^{-m}$. See for instance [\ref{FultHar}], Section~15.5.

\ssection
\ssectlabel{GIrreps}
Next we consider the reductive group $G = \Res_{K/k}\, H$ over~$k$. If $R$ is a commutative $k$-algebra, $G(R) = (B \otimes_k R)^*$. The set $\chargp(G)^+$ of dominant weights of~$G_\kbar$ is given by $\chargp(G)^+ = \oplus_{\sigma \in \Sigma(K)}\, \chargp(H)^+$. Via the identification $\chargp(H)^+ = \Lambda^+$ of~\ref{HIrreps}, we obtain an identification of~$\chargp(G)^+$ with the set 
$$
\chargp^+ = \bigoplus_{\sigma \in \Sigma(K)}\, \Lambda^+\, .
$$
The Galois group $\Gal(\kbar/k)$ acts on $\chargp^+ = \chargp(G)^+$ by its permutation of the summands; hence this action factors through an action of~$\Gamma$. By [\ref{Tits}], Thm.~7.2, the irreducible $k$-representations of~$G$ are indexed by the elements of $\chargp^+/\Gamma$. If $\xi$ is a $\Gamma$-orbit in~$\chargp^+$ we denote the corresponding irreducible representation of~$G$ by~$\rho_\xi$.

We have a natural isomorphism $G_\Ktilde \cong \prod_{\sigma \in \Sigma(K)}\, H_\sigma$, with $H_\sigma = H \otimes_{K,\sigma} \Ktilde$. The representation~$\rho_{\xi,\Ktilde}$ decomposes as a direct sum $\oplus_{\blambda \in \xi}\, \Psi_\blambda$, where $\Psi_\blambda$ is the external tensor product $\sqtensor_{\sigma \in \Sigma(K)}\, \psi_{\blambda(\sigma)}$. (Here $\blambda \in \chargp^+$ is viewed as a function $\Sigma(K) \to \Lambda^+$.)

Note that, since $G(k) = B^*$ is Zariski dense in~$G$, the representations~$\rho_\xi$, for $\xi \in \chargp^+/\Gamma$, are still irreducible and mutually non-equivalent as representations of the abstract group~$B^*$. 

\ssection
Choose a $k$-basis $\{\beta_1,\ldots,\beta_N\}$ for~$B$ (with $N=nd^2$). If $E$ is a commutative $k$-algebra, we call a map $r \colon B \to E$ a multiplicative homogeneous polynomial map over~$k$ of degree~$i$ if it has the following properties:
\item{(a)} $r$ is multiplicative, in the sense that $r(1) = 1$ and $r(b_1 b_2) = r(b_1) r(b_2)$ for all $b_1$, $b_2 \in B$;
\item{(b)} there exists a homogeneous polynomial $P \in E[t_1,\ldots,t_N]$ of degree~$i$ such that $r(c_1 \beta_1 + \cdots + c_N \beta_N) = P(c_1,\ldots,c_N)$ for all $c_1,\ldots,c_N \in k$.

\noindent
Note that the polynomial~$P$ in~(b) is uniquely determined, because $k$ is an infinite field.

Let $V$ be a finite dimensional $k$-vector space. Consider a multiplicative homogeneous polynomial map $r \colon B \to \End_k(V)$ over~$k$ of degree~$i$. If $R$ is a commutative $k$-algebra, define $r_R \colon B \otimes_k R \to \End_R(V \otimes_k R) = \End_k(V) 
\otimes_k R$ by the relation $r_R(c_1 \beta_1 + \cdots + c_N \beta_N) = P(c_1,\ldots,c_N)$, for $c_1,\ldots,c_N \in R$. Using that $r$ is multiplicative plus the fact that the field~$k$ is infinite, one easily shows that the map~$r_R$ is again multiplicative. Hence this construction defines an algebraic representation $\phi_r \colon G \to \GL(V)$ over~$k$. We refer to the representations of~$G$, or of $B^* = G(k)$, that are obtained in this manner as the polynomial representations of degree~$i$.

\ssection
\ssectlabel{XpolNotation}
Define a subset $\Lambda^\pol \subset \Lambda^+$ by the condition that $\lambda_d \geq 0$, i.e.,
$$
\Lambda^\pol = \bigl\{\lambda = (\lambda_1,\ldots,\lambda_d) \in \mZ^d \bigm| \lambda_1 \geq \lambda_2 \geq \cdots \geq \lambda_d \geq 0\bigr\}\, . \eqlabel{LambdaPolDef}
$$
Define $\chargp^\pol = \oplus_{\sigma \in \Sigma(K)}\, \Lambda^\pol$, which is a $\Gamma$-stable subset of~$\chargp^+$.

For $\blambda \in \chargp^\pol$, define $\wt{\blambda} =  \sum_{\sigma \in \Sigma(K)}\, \bigl|\blambda(\sigma)\bigr|$. As the map $\chargp^\pol \to \mZ_{\geq 0}$ given by $\blambda \mapsto \wt{\blambda}$ is $\Gamma$-invariant, it descends to a map $\wt{\ } \colon \chargp^\pol/\Gamma \to \mZ_{\geq 0}$.

\ssection
\ssectlabel{PolRepDec}
{\it Proposition. --- Let $\phi \colon B^*  \to \GL(V)$ be a polynomial representation of degree~$i$. Then there is a unique decomposition
$$
(V,\phi) = \bigoplus_{{\xi \in \chargp^{\pol}/\Gamma \atop \wt{\xi} = i}}\, (V^{(\xi)},\phi^{(\xi)}) \eqlabel{VphiDec}
$$
such that $(V^{(\xi)},\phi^{(\xi)})$ is isomorphic to a sum of copies of the irreducible representation~$\rho_\xi$.\/}
\medskip

\Proof
By construction, $\phi \colon B^*  \to \GL(V)$ is obtained from an algebraic representation $\phi_r \colon G \to \GL(V)$ by evaluation on $k$-rational points. The irreducible representations that occur in~$\phi_r$ are again polynomial of degree~$i$, and this property is preserved if we extend scalars to~$\Kbar$. Using the description of the representations~$\rho_{\xi,\Kbar}$ given in~\ref{HIrreps} and~\ref{GIrreps} we see that the only irreducible representations~$\rho_\xi$ that are polynomial of degree~$i$ are those with $\xi \in \chargp^\pol/\Gamma$ and $\wt{\xi} = i$. \QED

\ssection
\ssectlabel{RedNormExa}
{\it Example.\/} --- The reduced norm $\Nrd \colon B^* \to k^*$ is a polynomial representation of degree~$nd$. It corresponds to the $\Gamma$-orbit in~$\chargp^\pol$ that consists of the single element $\bnu \colon \Sigma(K) \to \Lambda^\pol$ with $\bnu(\sigma) = (1,\ldots,1)$ for all $\sigma \in \Sigma(K)$. If $\xi \in \chargp^\pol/\Gamma$ is the orbit of $\blambda \colon \Sigma(K) \to \Lambda^\pol$, the representation $\Nrd \otimes \rho_\xi$ is again polynomial; it corresponds to the $\Gamma$-orbit in~$\chargp^\pol$ of the sum $\bnu + \blambda$.

\ssection
\ssectlabel{PolRepDecInf}
{\it Remark.\/} --- We shall have to deal with multiplicative homogeneous polynomial maps $r \colon B \to \End_k(V)$ of degree~$i$ where $V$ is no longer assumed to have finite $k$-dimension. Then $V$ is the union of its finite dimensional subspaces~$V^\prime$ that are stable under all operators~$r(b)$ for $b \in B$; see for instance [\ref{DemGab}], \Romno2, \S 2, 3.1. Hence, we again have a decomposition~\eqref{VphiDec}, of course with the understanding that the $(V^{(\xi)},\phi^{(\xi)})$ will now in general be infinite sums of copies of~$\rho_\xi$. We refer to $V^{(\xi)}$ as the $\xi$-isotypic component of~$V$.

\section{Preliminaries on the action of endomorphisms on the Chow motive}{ChowAction}

\ssection
\ssectlabel{MotBasics}
Throughout this section, $F$ is a Dedekind ring and $S$ denotes a connected scheme that is smooth and of finite type over~$F$. Let $\Mot^0(S)$ be the category of Chow motives over~$S$ with respect to graded correspondences, as defined as in~[\ref{DenMur}],~1.6. Note that all results from~[\ref{DenMur}] are valid in the generality considered here; see [\ref{KunArChow}], Remark~1.1. (I thank K.~K\"unnemann for pointing this out to me.)

Let $\Var_S$ denote the category of smooth projective $S$-schemes. We have a contravariant functor $\Var_S \to \Mot^0(S)$, sending a smooth projective $X \to S$ to $R(X/S) = (X,[\transp \Gamma_\id],0)$ and sending $f \colon X \to Y$ to the morphism $[\transp \Gamma_f] \colon R(Y/S) \to R(X/S)$. The induced map on Chow rings is the pull-back map~$f^*$.

Let $X \to S$ be an abelian scheme of relative dimension~$g$ over~$S$. For $m \in \mZ$, let $\mult{m} \colon X \to X$ denote the multiplication by~$m$ map. By [\ref{DenMur}], Cor.~3.2 the relative motive $R(X/S)$ decomposes in $\Mot^0(S)$ as 
$$
R(X/S) = \bigoplus_{i=0}^{2g}\, R^i(X/S)\, , \eqlabel{DenMurDec}
$$
in such a way that $[\transp\Gamma_{\mult{m}}]$ acts on $R^i(X/S)$ as multiplication by~$m^i$. Let $[\Delta_{X/S}] = \sum_{i=0}^{2g} \pi_i$ be the corresponding decomposition of the diagonal of $X \times_S X$. If the context requires it we write $\pi_{X/S,i}$ instead of~$\pi_i$. It is convenient to define $R^i(X/S) = 0$ and $\pi_i=0$ if $i \notin \{0,\ldots,2g\}$.

If $f \colon X \to Y$ is a homomorphism of abelian schemes over~$S$ we have $[\transp\Gamma_f] \circ \pi_{Y/S,i} = \pi_{X/S,i} \circ [\transp\Gamma_f]$ for all~$i$.
\medskip

The first goal of this paper is to explain how, in the presence of non-trivial endomorphisms, the decomposition \eqref{DenMurDec} may be refined. As a first example we consider the case of a product of abelian schemes. The following result is an immediate consequence of the work of Deninger and Murre in~[\ref{DenMur}].

\ssection
\ssectlabel{DecProdAV}
{\it Proposition. --- Let $X_1,\ldots,X_r$ be abelian schemes over~$S$ with $X_\nu$ of relative dimension~$g_\nu$. Write $X = X_1 \times_S \cdots \times_S X_r$, let $g = g_1 + \cdots + g_r$ and 
$$
I_X = \bigl\{\bi = (i_1,\ldots,i_r) \in \mZ^r \bigm| 0 \leq i_\nu \leq 2g_\nu\bigr\}\, .
$$
For $\bm = (m_1,\ldots,m_r)$, let $\mult{\bm} \in \End(X/S)$ be given by $(x_1,\ldots,x_r) \mapsto (m_1x_1,\ldots,m_rx_r)$, and let $\bm^\bi = m_1^{i_1} \cdots m_r^{i_r}$. Then there is a unique decomposition 
$$
\bigl[\Delta_{X/S}\bigr] = \sum_{\bi \in I_X}\, \pi_\bi \eqlabel{DeltaDec}
$$
in $\End_{\Mot^0(S)}\bigl(R(X/S)\bigr) = \CH^g(X\times_S X)$ such that the elements $\pi_\bi$ are mutually orthogonal idempotents and such that $[\transp\Gamma_{\mult{\bm}}] \circ \pi_\bi = \bm^\bi \cdot \pi_\bi$ for all $\bm \in \mZ^r$ and $\bi \in I_X$. Moreover, $\pi_\bi \circ [\transp\Gamma_{\mult{\bm}}] = \bm^\bi \cdot \pi_\bi$ for all $\bm$ and~$\bi$. Corresponding to~{\rm \eqref{DeltaDec}} we have a decomposition
$$
R(X/S) = \bigoplus_{\bi\in I_X}\, R^\bi(X/S) 
$$
such that $[\transp\Gamma_{\mult{\bm}}]$ acts on $R^\bi(X/S)$ as multiplication by $\bm^\bi$.}
\medskip

\Proof
This follows from the main results of [\ref{DenMur}] by taking tensor products. We have $R(X/S) = R(X_1/S) \otimes \cdots \otimes R(X_r/S)$ in $\Mot^0(S)$. Now take $\pi_\bi = \pi_{X_1/S,i_1} \otimes \cdots \otimes \pi_{X_r/S,i_r}$ for $\bi = (i_1,\ldots,i_r) \in I_X$. \QED

\ssection
\ssectlabel{MotFourDual}
{\it Example.\/} --- (Cf.\ [\ref{KunChow}], (3.1.2)(\romno2).) Let $X$ and~$Y$ be abelian schemes over~$S$ with $X$ of relative dimension~$g$. If $z \in \CH(X \times_S Y)$ we have a decomposition $z = \sum z_{i,j}$ such that $[m,n]^*(z_{i,j}) = m^i n^j\cdot z_{i,j}$ for all integers $m$ and~$n$. It follows from the relations in [\ref{DenMur}], Prop.~1.2.1, together with the motivic Poincar\'e duality $\transp{\pi_i} = \pi_{2g-i}$ that $z_{i,j} = \pi_{Y/S,j} \circ z \circ \pi_{X/S,2g-i}$.

We apply this with $Y = X^\avdual$, the dual of~$X$. Let $\ell = \ell_X \in \CH^1(X \times_S X^\avdual)$ be the first Chern class of the Poincar\'e bundle, and recall that we use the divided power notation $\ell^{[n]} = \ell^n/n!$. Then $\ell = \ell_{1,1}$; hence, $\ell^{[i]} = \pi_i(X^\avdual/S) \circ \ell^{[i]} \circ \pi_{2g-i}(X/S)$. Now use the Mukai-Beauville relation $\Four^\avdual \circ \Four = (-1)^g [-1]^*$ and view $\ell^{[i]} \in \CH^i(X\times_S X^\avdual)$ as a morphism from $R(X/S) = \oplus R^j(X/S)$ to $R(X^\avdual/S)\bigl(g-i\bigr) = \oplus R^j(X^\avdual/S)\bigl(i-g\bigr)$. It follows that the only non-zero component of this morphism is a Fourier duality isomorphism
$$
\ell^{[i]} \colon R^{2g-i}(X/S) \isomarrow R^i(X^\avdual/S)\bigl(i-g\bigr)\, . \eqlabel{lnovern!}
$$
(The interpretation is that the dual abelian scheme is the Poincar\'e dual of~$X$. Indeed, combining~\eqref{lnovern!} with the motivic Poincar\'e duality $R^i(X/S)^\vee = R^{2g-i}(X/S)\bigl(g\bigr)$ we find that $R^i(X^\avdual/S) \cong R^i(X/S)^\vee\bigl(-i\bigr)$.)

\ssection
\ssectlabel{AVwithEnd}
With $S$ as in~\ref{MotBasics}, consider an abelian scheme $X \to S$ of relative dimension $g>0$. If $m$ is a nonzero integer, it follows from the results discussed above that the endomorphisms $[\transp\Gamma_{\mult{m}}]$ and $[\Gamma_{\mult{m}}]$ of $R(X/S)$ are invertible. Every element~$\alpha$ of the endomorphism algebra $D = \End^0(X/S)$ can be written in the form $\alpha = f/m$ for some $f \in \End(X/S)$ and some integer $m \neq 0$. We then define classes $[\transp\Gamma_\alpha]$ and $[\Gamma_\alpha]$ in $\CH^g(X\times_S X)$ by 
$$
[\transp\Gamma_\alpha] = [\transp\Gamma_f] \circ [\transp\Gamma_{\mult{m}}]^{-1}
\; ,\qquad
[\Gamma_\alpha] = [\Gamma_f] \circ [\Gamma_{\mult{m}}]^{-1}\, .
$$
In what follows we denote the induced endomorphisms of~$\CH(X)$ simply by $\alpha^*$ and~$\alpha_*$.

For $i\geq 0$ we have $\pi_i \circ [\transp\Gamma_\alpha] = [\transp\Gamma_\alpha] \circ \pi_i$. Define $r^{(i)} \colon D^\opp \to \End_{\Mot^0(S)}\bigl(R^i(X/S)\bigr)$ by $r^{(i)}(\alpha) = [\transp\Gamma_\alpha] \circ \pi_i$. The maps~$r^{(i)}$ are multiplicative but not, in general, additive. In particular, the group $D^{\opp,*}$ acts on~$R^i(X/S)$ by automorphisms.

\ssection
\ssectlabel{risHomom}
{\it Proposition. --- The map $r^{(1)} \colon D^\opp \to \End_{\Mot^0(S)}\bigl(R^1(X/S)\bigr)$ is an isomorphism of $\mQ$-algebras.}
\medskip

\Proof
See [\ref{Kings}], Prop.~2.2.1. 
\QED

\ssection
\ssectlabel{riismhp}
{\it Corollary. --- The map $r^{(i)} \colon D^\opp \to \End_{\Mot^0(S)}\bigl(R^i(X/S)\bigr)$ defined in\/~{\rm \ref{AVwithEnd}} is a multiplicative homogeneous polynomial map over~$\mQ$ of degree~$i$.}
\medskip

\Proof
We already know that $r^{(i)}$ is multiplicative. Taking the isomorphism $R^i(X/S) \isomarrow \wedge^i R^1(X/S)$ of [\ref{KunChow}], Thm.~(3.3.1), as an identification, the map $r^{(i)}$ is the composition of the homomorphism~$r^{(1)}$ with the map $\End_{\Mot^0(S)}\bigl(R^1(X/S)\bigr) \to \End_{\Mot^0(S)}\bigl(R^i(X/S)\bigr)$ that sends an endomorphism~$h$ of~$R^1(X/S)$ to the induced endomorphism $\wedge^i h = h\wedge \cdots \wedge h$ of~$R^i(X/S)$. It follows that $r^{(i)}$ is a homogeneous polynomial map of degree~$i$. \QED

\section{Duality}{Duality}

\ssection
\ssectlabel{ThreeBs}
Again let $X \to S$ be an abelian scheme of relative dimension~$g$. We assume that the endomorphism algebra $D = \End^0(X/S)$ is a simple $\mQ$-algebra of finite dimension. (For the general case see~\ref{GeneralCase}.) Let $K$ be the center of~$D$. Let $n = [K:\mQ]$ and $d = \dim_K(D)^{1/2}$. Let $\Sigma(K)$ be the set of ring homomorphisms $K \to \Qbar$, let $\Ktilde \subset \Qbar$ denote the normal closure of~$K$ inside~$\Qbar$, and write $\Gamma = \Gal(\Ktilde/\mQ)$.

We apply the theory of Section~\ref{RepTheory} with $k=\mQ$ and three different choices for~$B$, to be discussed in more detail below. In each case $B$~is central simple of dimension~$d^2$ over~$K$. The meaning of $\Sigma(K)$ and~$\Gamma$ is the same in all cases and the notation we use is consistent with the notation introduced in Section~\ref{RepTheory}. In each case we index the irreducible algebraic representations of~$B^*$ by $\chargp^+/\Gamma$, following the method discussed in \ref{HRootDatum}--\ref{GIrreps}.

Let us now give some more details about the group actions we consider.

\item{(a)} We shall mostly take the cohomological perspective. In this case we consider $B = D^\opp$, acting on~$\CH(X)$ through the operators~$f^*$. Let $H$ denote the reductive group over~$K$ with $H(R) = (D^\opp \otimes_K R)^*$ and let $G = \Res_{K/\mQ}\, H$. For $\lambda \in \Lambda^+$, let $\psi_\lambda$ be the corresponding irreducible representation of~$H$ over~$K$. For $\xi \in \chargp^+/\Gamma$, let $\rho_\xi$ be the corresponding irreducible representation of $G(\mQ) = D^{\opp,*}$ over~$\mQ$.

\item{(b)} In order to describe Poincar\'e duality we need the homological perspective, letting $B=D$ act on~$\CH(X)$ through the operators~$f_*$. Let $H^\prime$ be the reductive group over~$K$ with $H^\prime(R) = (D \otimes_K R)^*$ and let $G^\prime = \Res_{K/\mQ}\, H^\prime$, which is the opposite of the group~$G$. For $\lambda \in \Lambda^+$, let $\psi^\prime_\lambda$ be the corresponding irreducible representation of~$H^\prime$ over~$K$. For $\xi \in \chargp^+/\Gamma$, the corresponding irreducible representation of $G^\prime(\mQ) = D^*$ over~$\mQ$ is denoted by~$\rho^\prime_\xi$.

\item{(c)} Let $X^\avdual \to S$ be the dual abelian scheme and let $D^\avdual = \End(X^\avdual/S) \otimes \mQ$. If $f$ is an endomorphism of $X/S$, let $f^\avdual \colon X^\avdual \to X^\avdual$ denote the dual endomorphism. The map $f \mapsto f^\avdual$ gives an isomorphism of $\mQ$-algebras $D \isomarrow D^{\avdual,\opp}$ and we use this to identify the center of~$D^{\avdual,\opp}$ with~$K$. (This may lead to confusion; see~\ref{KCaution}.) For the rest the pattern is the same as in~(a). We consider $\CH(X^\avdual)$ as a representation of $D^{\avdual,\opp,*}$, with $g \in D^{\avdual,\opp}$ acting as~$g^*$. For $\xi \in \chargp^+/\Gamma$, let $\rho^\avdual_\xi$ be the corresponding irreducible representation of $D^{\avdual,\opp,*}$ over~$\mQ$.

\ssection
\ssectlabel{HWeightLem}
{\it Lemma. --- Let $\lambda = (\lambda_1,\lambda_2,\ldots,\lambda_d)$ be an element of~$\Lambda^+$. Then the representation~$\tau$ of~$H^\prime$ over~$K$ given by $\tau(h) = \psi_\lambda(h^{-1})$ is isomorphic to~$\psi^\prime_\mu$, where $\mu = (-\lambda_d,\ldots,-\lambda_1)$.}
\medskip

\Proof
It is clear that $\tau$ is an irreducible representation of~$H^\prime$.
As the representations are determined by their highest weights, we may work over~$\Kbar$. Choose an isomorphism of $\Kbar$-algebras $a \colon D_\Kbar^\opp \isomarrow M_d(\Kbar)$, and define $a^\prime \colon D_\Kbar \isomarrow M_d(\Kbar)$ by $a^\prime(\delta) = \transp{a(\delta)}$, the transpose of $a(\delta)$. Let $\xi \colon H_\Kbar \isomarrow \GL_{d,\Kbar}$ and $\xi^\prime \colon H^\prime_\Kbar \isomarrow \GL_{d,\Kbar}$ be the induced isomorphisms of algebraic groups. Via these isomorphisms we can view both~$\psi_\lambda$ and $\tau$ as representations of~$\GL_{d,\Kbar}$; in other words, we consider $\psi_\lambda \circ \xi^{-1}$ and $\tau \circ (\xi^\prime)^{-1}$. In both cases the highest weight is taken with regard to the diagonal torus~$T$ and the upper triangular Borel $Q \subset \GL_d$. We have
$$
\bigl(\tau \circ (\xi^\prime)^{-1}\bigr)(g) = \bigl(\psi_\lambda \circ \xi^{-1}\bigr)(\transp{g}^{-1})\, .
$$
Let $\beta$ be the automorphism of~$\GL_d$ given by $g \mapsto \transp{g}^{-1}$. Then $\beta(T) = T$ and $\beta(Q) = Q^-$, the lower triangular Borel subgroup. If $A \in \GL_d(K)$ is the anti-diagonal matrix with all anti-diagonal coefficients equal to~$1$, the inner automorphism $\Inn(A)$ transforms $(T,Q^-)$ back to~$(T,Q)$, and the effect of $\Inn(A) \circ \beta$ on the character group of~$T$ is given by $e_i \mapsto -e_{d-i}$. Hence if $\psi_\lambda \circ \xi^{-1}$ has highest weight $\lambda_1 e_1 + \cdots + \lambda_d e_d$, the highest weight of $\tau \circ (\xi^\prime)^{-1}$ is $-\lambda_d e_1 - \cdots - \lambda_1 e_d$.
\QED

\ssection
\ssectlabel{xiwdualDef}
{\it Notation.\/} --- For $\lambda = (\lambda_1,\ldots,\lambda_d)$ in~$\Lambda^+$ define 
$$
\lambda^\wdual = \left({2g\over nd}-\lambda_d,\ldots,{2g\over nd}-\lambda_1\right)\, .
$$
Note that $2g/nd$ is an integer; see [\ref{MAV}], Chap.~19, Corollary to Thm.~4. Hence $\lambda^\wdual$ is again an element of~$\Lambda^+$. For $\blambda \in \chargp^+$, define $\blambda^\wdual \in \chargp^+$ by the rule $\blambda^\wdual(\sigma) = \blambda(\sigma)^\wdual$. For $\xi \in \chargp^+/\Gamma$, let $\xi^\wdual$ denote the $\Gamma$-orbit consisting of the elements~$\blambda^\wdual$, for $\blambda \in \xi$. Note that $\wt{\xi^\wdual} = 2g-\wt{\xi}$.

\ssection
\ssectlabel{FDualEffect}
{\it Proposition. --- Let $V \subset \CH(X)$ be an irreducible subrepresentation of~$D^{\opp,*}$ that is isomorphic to~$\rho_\xi$.

{\rm (\romno1)} The subspace $V\subset \CH(X)$ is stable under the action of the operators~$f_*$, for $f \in D$, and $V$ is isomorphic to $\rho^\prime_{\xi^\wdual}$ as a representation of~$D^*$. 

{\rm (\romno2)} Let $\Four \colon \CH(X) \isomarrow \CH(X^\avdual)$ be the Fourier transform. Then $\Four(V) \subset \CH(X^\avdual)$ is an irreducible subrepresentation of~$D^{\avdual,\opp,*}$ that is isomorphic to~$\rho^\avdual_{\xi^\wdual}$.}
\medskip

\Proof
(\romno1) Let $f \in D^{\opp,*}$. Then $f$ is a quasi-isogeny of~$X$ to itself. Its degree $\deg(f)$ equals $\Nrd(f)^{(2g/nd)}$, where $\Nrd\colon D^{\opp,*} \to \mQ^*$ is the reduced norm character. (See~\ref{RedNormExa}.) For $z \in \CH(X)$ we have the relation $f_*(z) = \deg(f) \cdot (1/f)^*\bigl(z\bigr)$. Now use \ref{RedNormExa} and~Lemma~\ref{HWeightLem}.

(\romno2) For $f \in D$ and $z \in \CH(X)$ we have the relation $\Four\bigl(f_*(z)\bigr) = f^{\avdual,*}\bigl(\Four(z)\bigr)$. So (\romno2) follows from~(\romno1).
\QED

\ssection
\ssectlabel{KCaution} 
{\it Caution.\/} --- The field~$K$ is either totally real or a CM field. In (\romno2) of the Proposition, it is important that we identify $K$ with the center of~$D^{\avdual,\opp}$ via the isomorphism $D \isomarrow D^{\avdual,\opp}$ given by $f \mapsto f^\avdual$. If we choose a polarization $\theta\colon X \to X^\avdual$, the resulting isomorphism $D \isomarrow D^\avdual$ gives the complex conjugate identification of~$K$ with the center of~$D^{\avdual,\opp}$. Under that identification, the Fourier dual of a $D^{\opp,*}$-subrepresentation $V \subset \CH(X)$ of type~$\rho_\xi$ is a $D^{\avdual,\opp,*}$-subrepresentation $\Four(V) \subset \CH(X^\avdual)$ of type $\rho^\avdual_{\bar{\xi}^\wdual}$, where $\bar{\xi}^\wdual \in \chargp^\adm/\Gamma$ is the complex conjugate of~$\xi^\wdual$.

\section{Motivic decomposition}{MotivicDec}
\medskip

\noindent
We retain the notation and assumptions of~\ref{ThreeBs}.

\ssection
\ssectlabel{LadmDef}
Define a subset $\Lambda^\adm \subset \Lambda^\pol$ of ``admissible'' elements by the condition that $(2g/nd) \geq \lambda_1$; so,
$$
\Lambda^\adm = \Bigl\{\lambda = (\lambda_1,\ldots,\lambda_d) \in \mZ^d \Bigm| {2g\over nd} \geq \lambda_1 \geq \lambda_2 \geq \cdots \geq \lambda_d \geq 0 \Bigr\}\, .
$$
Define $\chargp^\adm = \oplus_{\sigma \in \Sigma(K)}\, \Lambda^\adm$, which is a $\Gamma$-stable subset of~$\chargp^+$. Note that $0 \leq \wt{\xi} \leq 2g$ for all $\xi \in \chargp^\adm/\Gamma$. If $\blambda \in \chargp^\adm$ then $\blambda^\wdual$ is an element of~$\chargp^\adm$, too; hence $\xi \mapsto \xi^\wdual$ is an involutive automorphism of~$\chargp^\adm/\Gamma$. 

\ssection
\ssectlabel{CHDecThm}
{\it Theorem. --- There is a unique decomposition
$$
\CH(X) = \bigoplus_{\xi \in \chargp^\adm/\Gamma}\, \CH^{(\xi)}(X) \eqlabel{CHRiDec}
$$
as a representation of~$D^{\opp,*}$, such that $\CH^{(\xi)}(X)$ is $\xi$-isotypic. The subspace $\CH\bigl(R^i(X/S)\bigr) \subset \CH(X)$ is the direct sum of the $\CH^{(\xi)}(X)$ with $\wt{\xi}=i$. For $\xi \in \chargp^\adm/\Gamma$, the Fourier transform~$\cF$ restricts to an isomorphism $\cF \colon \CH^{(\xi)}(X) \isomarrow \CH^{(\xi^\wdual)}(X^\avdual)$.}
\medskip

\Proof
By \ref{riismhp} we can apply Prop.~\ref{PolRepDec}. This gives a decomposition of~$\CH\bigl(R^i(X/S)\bigr)$ as a direct sum of subspaces $\CH^{(\xi)}\bigl(R^i(X/S)\bigr)$ for $\xi \in \chargp^\pol/\Gamma$ with $\wt{\xi} = i$. (Cf.~\ref{PolRepDecInf}.) If $\CH^{(\xi)}\bigl(R^i(X/S)\bigr) \neq 0$ then it follows from Prop.~\ref{FDualEffect} that $\xi^\wdual$ lies in the subset $\chargp^\pol/\Gamma \subset \chargp^+/\Gamma$. This implies that $\xi \in \chargp^\adm/\Gamma$. The last assertion is immediate from Prop.~\ref{FDualEffect}(\romno2). \QED
\medskip

As the decomposition of the Chow group has good functorial properties, we obtain from it a motivic decomposition.

\ssection
\ssectlabel{MotDecRefined}
{\it Theorem. ---  {\rm (\romno1)} There is a unique decomposition
$$
R(X/S) = \bigoplus_{\xi \in \chargp^\adm/\Gamma}\, R^{(\xi)}(X/S)\, , \eqlabel{RX/SDec}
$$
in $\Mot^0(S)$ that is stable under the action of~$D^{\opp,*}$ and has the property that for any $M$ in~$\Mot^0(S)$ the $D^{\opp,*}$-representation $\Hom_{\Mot^0(S)}\bigl(M,R^{(\xi)}(X/S)\bigr)$ is $\xi$-isotypic. The submotive $R^i(X/S)$ is the direct sum of the $R^{(\xi)}(X/S)$ with $\wt{\xi} = i$.

{\rm (\romno2)} For $\xi \in \chargp^\adm/\Gamma$ the subspace $\CH\bigl(R^{(\xi)}(X/S)\bigr) \subset \CH(X)$ is the $\xi$-isotypic component $\CH^{(\xi)}(X) \subset \CH(X)$ of\/~{\rm \eqref{CHRiDec}}.

{\rm (\romno3)} Let $\delta_\xi$ be the idempotent in  $\CH^g(X \times_S X) = \End_{\Mot^0(S)}\bigl(R(X/S)\bigr)$ that defines the submotive~$R^{(\xi)}(X/S)$, so that $[\Delta_{X/S}] = \sum_{\xi \in \chargp^\adm/\Gamma}\, \delta_\xi$ is the decomposition of the diagonal that corresponds with\/~{\rm \eqref{RX/SDec}}. Then $\transp{\delta_\xi} = \delta_{\xi^\wdual}$; hence 
$$
R^{(\xi)}(X/S)^\vee = R^{(\xi^\wdual)}(X/S)\bigl(g\bigr)\, .
$$

{\rm (\romno4)} The motivic Fourier duality $R^{2g-i}(X/S) \isomarrow R^i(X^\avdual/S)\bigl(i-g\bigr)$ of\/~{\rm \eqref{lnovern!}} is the direct sum of isomorphisms
$$
R^{(\xi)}(X/S) \isomarrow R^{(\xi^\wdual)}(X^\avdual/S)\bigl(i-g\bigr)
$$
for $\xi \in \chargp^\adm/\Gamma$ with $\wt{\xi} = 2g-i$.}
\medskip

\Proof
Let $M=(Y,p,m)$ be an object of~$\Mot^0(S)$. We first assume $Y$ is connected. Viewing $Y \times_S X$ as an abelian scheme over~$Y$ via the first projection, Thm.~\ref{CHDecThm} gives us an isotypic decomposition $\CH(Y\times_S X) = \oplus_{\xi \in \chargp^\adm/\Gamma}\, \CH^{(\xi)}(Y \times_S X)$. Note that $\alpha \in D$ acts on $\CH(Y\times_S X)$ as $(\id_Y \times \alpha)^*$. 

By definition, $\Hom_{\Mot^0(S)}\bigl(M,R(X/S)\bigr) = \CH^{d-m}(Y\times_S X) \circ p$, with $d=\dim(Y/S)$. By [\ref{DenMur}], Prop.~1.2.1, we have $(1 \times \alpha)^*\bigl(\zeta \circ p\bigr) = (1 \times \alpha)^*\zeta \circ p$ for all $\alpha \in D$ and $\zeta \in \CH(Y\times_S X)$; this just says that the endomorphism of $\CH(Y\times_S X)$ given by $\zeta \mapsto \zeta \circ p$ is $D^{\opp,*}$-equivariant. Hence we have a decomposition
$$
\Hom_{\Mot^0(S)}\bigl(M,R(X/S)\bigr) = \bigoplus_{\xi \in \chargp^\adm/\Gamma}\, H^{(\xi)} \eqlabel{HomMXDec}
$$
with $H^{(\xi)}$ a $\xi$-isotypic representation of~$D^{\opp,*}$. If $Y$ is not connected we obtain the same conclusion by first decomposing~$M$ as a finite sum of motives on connected $S$-schemes. Moreover, it is straightforward to verify that the decomposition~\eqref{HomMXDec} is functorial. By Yoneda, together with the fact that $\Mot^0(S)$ is pseudo-abelian, we therefore have a decomposition~\eqref{RX/SDec} with the properties stated in~(\romno1). 

Part~(\romno2) of the theorem follows from~(\romno1) by taking $M = \unitmot(-j)$ for various~$j$, and (\romno4) follows from the last assertion of Thm.~\ref{CHDecThm}, again using a Yoneda argument. 

For (\romno3) we first recall from \ref{ThreeBs}(c) that we have a natural isomorphism $\tau \colon D^* \cong D^{\avdual,\opp,*}$. On $R^i(X/S)^\vee$ we have an action of~$D^*$. On $R^i(X^\avdual/S)\bigl(i\bigr)$ we have an action of~$D^{\avdual,\opp,*}$. Further, the isomorphism $R^i(X/S)^\vee \isomarrow R^i(X^\avdual/S)\bigl(i\bigr)$ of~\ref{MotFourDual} is equivariant with respect to~$\tau$. (Cf.\ the proof of \ref{FDualEffect}(\romno2).) With these remarks, (\romno3) follows from~(\romno4).
\QED

\ssection
{\it Remark.\/} --- Viewing $X \times_S X$ as an abelian scheme over~$X$ via the first projection we have a decomposition $\CH(X\times_S X) = \oplus_{\xi \in \chargp^\adm/\Gamma}\, \CH^{(\xi)}(X \times_S X)$. The projectors~$\delta_\xi$ that give~\eqref{RX/SDec} are the components of the identity~$[\transp\Gamma_\id]$ in this decomposition.

\ssection
\ssectlabel{DquatExa}
{\it Example.\/} --- Suppose $D$ is a quaternion algebra with center~$\mQ$. In this case $\chargp(G)^\adm/\Gamma$ is the set of pairs $\lambda = (\lambda_1,\lambda_2)$ with $g \geq \lambda_1 \geq \lambda_2 \geq 0$. Viewing~$D^{\opp,*}$ as an inner form of~$\GL_2$ over~$\mQ$, the irreducible representation~$\rho_\lambda$ associated with~$\lambda$ (which in this case is the same as the representation~$\psi_\lambda$ of~\ref{HIrreps}) is a $\mQ$-form of $d(\lambda)$ copies of the representation $\Sym^{\lambda_1-\lambda_2}(V) \otimes \det^{\otimes \lambda_2}$, where $V$ is the standard representation of~$\GL_2$ and where
$$
d(\lambda) = \cases{1 & if $\lambda_1-\lambda_2$ is even;\cr 2 & if $\lambda_1-\lambda_2$ is odd.} \eqlabel{dlambdaDef}
$$

For $0 \leq i \leq g$ we obtain a decomposition
$$
R^i(X/S) = R^{(i,0)} \oplus R^{(i-1,1)} \oplus \cdots \oplus R^{(\nu,i-\nu)} \quad \hbox{with $\nu = \lfloor i/2 \rfloor$.}
$$
For $g \leq i \leq 2g$ the decomposition takes the form
$$
R^i(X/S) = R^{(g,i-g)} \oplus R^{(g-1,i+1-g)} \oplus \cdots \oplus R^{(g-\nu,i+\nu-g)} \quad \hbox{with $\nu = \lfloor (2g-i)/2 \rfloor$.}
$$
Fourier duality exchanges $R^{(\lambda_1,\lambda_2)}(X/S)$ and $R^{(g-\lambda_2,g-\lambda_1)}(X^\avdual/S)$. By looking at cohomology we can see that in general all summands $R^{(\lambda_1,\lambda_2)}$ in the indicated range are non-zero.
\medskip

A formal consequence of Thm.~\ref{MotDecRefined} is that we also get a decomposition of any realization of $R(X/S)$ on a $\mQ$-vector space. 

\ssection
\ssectlabel{Dopp*Realis}
{\it Corollary. --- Let $\Vect_\mQ$ be the category of $\mQ$-vector spaces. If $\Phi \colon \Mot^0(S) \to \Vect_\mQ$ is a $\mQ$-linear functor, $\Phi\bigl(R(X/S)\bigr) = \oplus_{\xi \in \chargp^\adm/\Gamma} \, \Phi\bigl(R^{(\xi)}(X/S)\bigr)$ and $\Phi\bigl(R^{(\xi)}(X/S)\bigr)$ is $\xi$-isotypic as a representation of~$D^{\opp,*}$.}
\medskip

\Proof
Write $R^{(\xi)} = R^{(\xi)}(X/S)$. The endomorphism algebra $\End(R^{(\xi)})$ is a representation of~$D^{\opp,*}$ through its action on the target, and by Thm.~\refn{MotDecRefined}(\romno1) this representation is $\xi$-isotypic. Let $E_\xi \subset \End(R^{(\xi)})$ be the $D^{\opp,*}$-subrepresentation generated by the identity. Equivalently, $E_\xi$ is the image of the group algebra $\mQ[D^{\opp,*}]$ in $\End(R^{(\xi)})$. If $u \in \Phi(R^{(\xi)})$, the $D^{\opp,*}$-submodule of $\Phi(R^{(\xi)})$ generated by~$u$ is a quotient of~$E_\xi$; hence this submodule is again $\xi$-isotypic. 
\QED

\ssection
{\it Example.\/} --- For the higher Chow groups (with $\mQ$-coefficients) we have 
$$
\CH(X;j) = \bigoplus_{\xi \in \chargp^\adm/\Gamma}\, \CH\bigl(R^{(\xi)}(X/S);j\bigr)
$$
and $\CH\bigl(R^{(\xi)}(X/S);j\bigr)$ is $\xi$-isotypic as a representation of~$D^{\opp,*}$. 

Depending on the context we can draw similar conclusions for cohomology. For instance, if the ground field~$F$ is~$\mC$ and if $q\colon X \to S$ is the structural morphism, the variation of Hodge structure $\mV = R^n q_*\mQ_X$ decomposes as a direct sum $\oplus_{\xi \in \chargp^\adm/\Gamma} \mV^{(\xi)}$ where $\mV^{(\xi)} \subset \mV$ is cut out by the projector~$\delta_\xi$ and is $\xi$-isotypic as a sheaf of $D^{\opp,*}$-modules.

If we have a cohomology theory with coefficients in a field~$\mF$ of characteristic~$0$, we can in general only conclude that the cohomology of $R^{(\xi)}(X/S)$ is a quotient of a sum of copies of~$\rho_{\xi,\mF}$. For instance, if $E$ is a supersingular elliptic curve over~$\ol{\mF}_p$, in which case $D$ is a quaternion algebra over~$\mQ$, there is a unique class $\xi \in \chargp^\adm/\Gamma$ with $\wt{\xi} = 1$ (see~\ref{DquatExa}) and $\rho_\xi$ has dimension~$4$; so the $\ell$-adic cohomology $H^1(E,\mQ_\ell)$ is only ``half'' a copy of~$\rho_{\xi,\mQ_\ell}$.

\ssection
\ssectlabel{GeneralCase}
{\it Remark.\/} --- So far we have assumed that the endomorphism algebra $\End^0(X/S)$ is a simple $\mQ$-algebra. This assumption was made for simplicity of exposition and is not essential. 

Let $\eta$ be the generic point of~$S$. Then $\End^0(X/S) = \End^0(X_\eta/\eta)$, so the assumption that $\End^0(X/S)$ is simple just means that the generic fiber~$X_\eta$ is isogenous to a power of a simple abelian variety over~$\eta$. We do not know if this implies that~$X$, as an abelian scheme over~$S$, is isogenous to a power of a simple abelian scheme. Similarly, if, in the general case, $\End^0(X/S) = D_1 \times \cdots \times D_r$ is the decomposition of the endomorphism algebra as a product of simple factors, we do not know if $X$ is necessarily isogenous to a product $Y_1 \times \cdots \times Y_r$ with $\End^0(Y_j) = D_j$. Nonetheless, the motive $R^1(X/S)$ behaves as if this were the case, since by Lemma~\ref{risHomom} we have a decomposition
$$
R^1(X/S) = M_1 \oplus \cdots \oplus M_r
$$
in $\Mot^0(S)$ with $D_j$ acting trivially on the factors~$M_i$ for $i \neq j$. All further arguments go through without essential changes, and Thms.~\ref{CHDecThm} and~\ref{MotDecRefined} are valid in this more general context, except that we have to replace the index set $\chargp^\adm/\Gamma$ by the product $\prod_{j=1}^r\, \chargp_j^\adm/\Gamma_j$ of the index sets associated to the factors~$D_j$. We leave it to the reader to write out the details.

Let us also note that, instead of taking $D = \End^0(X/S)$, we may work with a semisimple subalgebra $D \subset \End^0(X/S)$. In fact, taking a smaller algebra may give a finer motivic decomposition. For example, suppose $X = Y \times Y$ for some abelian scheme $Y/S$ with $\End(Y/S) = \mZ$. Then the decomposition of~\ref{DecProdAV} is finer than the decomposition of~$R(X/S)$ we obtain by applying~\ref{MotDecRefined} to~$X$, taking $D = \End^0(X/S)$. However, the finer decomposition in~\ref{DecProdAV} does not give information on how $\GL_2(\mQ)$ acts; it only takes into account the action of the diagonal subgroup $\mQ^* \times \mQ^*$.

\ssection
\ssectlabel{MotDecKtilde}
{\it Remark.\/} --- There is another, perhaps more elementary, way to obtain a motivic decomposition of $R(X/S)$, which coincides with \eqref{RX/SDec} if $D=K$ but which in general is coarser. For this we need to work in the category $\Mot^0(S;\Ktilde)$ of relative Chow motives with coefficients in the normal closure~$\Ktilde$. Write $R^i(X/S;\Ktilde)$ for the image of $R^i(X/S)$ under the natural functor $\Mot^0(S) \to \Mot^0(S;\Ktilde)$.

Let $D_\Ktilde = D \otimes_\mQ \Ktilde$. Then $D_\Ktilde = \prod_{\sigma\in \Sigma(K)}\, D_\sigma$, where $D_\sigma = D \otimes_{K,\sigma} \Ktilde$. Let $1 = \sum e_\sigma$ be the corresponding decomposition of $1 \in D_\Ktilde$ as a sum of idempotents. By~Prop.~\ref{risHomom} we have an algebra homomorphism $r_\Ktilde\colon D_\Ktilde^\opp \to \End_{\Mot^0(S;\Ktilde)}\bigl(R^1(X/S;\Ktilde)\bigr)$. This gives a decomposition $R^1(X/S;\Ktilde) = \oplus_{\sigma \in \Sigma(K)}\, R_\sigma$, where $R_\sigma$ is the submotive of $R^1(X/S;\Ktilde)$ cut out by the idempotent~$r_\Ktilde(e_\sigma)$. 

Let $\bJ = (\mZ_{\geq 0})^{\Sigma(K)}$, and for $i \geq 0$ define a subset $\bJ(i) \subset \bJ$ by
$$
\bJ(i) = \bigl\{\bj \colon \Sigma(K) \to \mZ_{\geq 0} \bigm| |\bj| = i \bigr\}\, ,
$$
where $|\bj| = \sum_{\sigma \in \Sigma(K)}\, \bj(\sigma)$. Taking exterior powers and using K\"unnemann's isomorphism $\wedge^i R^1(X/S) \isomarrow R^i(X/S)$, we obtain decompositions
$$
R^i(X/S;\Ktilde) = \bigoplus_{\bj \in \bJ(i)}\, R^{\{\bj\}}(X/S;\Ktilde)
\qquad\hbox{such that}\quad
R^{\{\bj\}}(X/S;\Ktilde) \cong \bigotimes_{\sigma \in \Sigma(K)}\, \Bigl(\wedge^{\bj(\sigma)}\, R_\sigma\Bigr)\, .
$$
(The calculation of the exterior powers works as expected; cf.\ [\ref{DelCatTens}], Section~1.) Fixing $i \geq 0$, let $1 = \sum_{\bj \in \bJ(i)}\, \tilde{e}_\bj$ be the corresponding decomposition of $1 \in \End_{\Mot^0(S;\Ktilde)}\bigl(R^i(X/S;\Ktilde)\bigr)$ as a sum of idempotents. The Galois group~$\Gamma$ acts on~$\bJ(i)$ and on the endomorphism algebra of the motive $R^i(X/S;\Ktilde)$. If $\gamma \in \Gamma$ sends $\bj \in \bJ(i)$ to~$\bj^\prime$ then ${}^\gamma \tilde{e}_\bj = \tilde{e}_{\bj^\prime}$. Hence if $\eta$ is a $\Gamma$-orbit in~$\bJ(i)$, the sum $\sum_{\bj \in \eta}\, \tilde{e}_\bj$ is an idempotent in $\End_{\Mot^0(S)}\bigl(R^i(X/S)\bigr)$. This gives us a decomposition
$$
R^i(X/S) = \bigoplus_{\eta \in \bJ(i)/\Gamma}\, R^{\{\eta\}}(X/S)
$$
in $\Mot^0(S)$ such that $R^{\{\eta\}}(X/S;\Ktilde) = \oplus_{\bj \in \eta}\, R^{\{\bj\}}(X/S;\Ktilde)$. 

To describe the relation with the decomposition in~\eqref{RX/SDec}, consider the map $v\colon \chargp^\adm/\Gamma \to \bJ/\Gamma$ that sends the $\Gamma$-orbit of $\blambda \in \chargp^\adm$ to the $\Gamma$-orbit of the function $\sigma \mapsto |\blambda(\sigma)|$. By analyzing how the groups $D_\sigma^{\opp,*}$ act, we find that $R^{\{\eta\}}(X/S) = \oplus\, R^{(\xi)}(X/S)$, where the sum runs over the classes $\xi \in \chargp^\adm/\Gamma$ such that $v(\xi) = \eta$.

\section{A Lie algebra action on the motive}{LieAction}
\medskip

\noindent
As before we consider an abelian scheme $X/S$. We no longer assume that the endomorphism algebra $D = \End^0(X/S)$ is simple. Write $\gl(X)$ for $D$ as a Lie algebra and define a Lie subalgebra $\sp(X\times_S X^\avdual) \subset \gl(X\times_S X^\avdual)$ by
$$
\sp(X\times_S X^\avdual) = \Biggl\{ 
\pmatrix{\alpha & \beta\cr \gamma & -\alpha^\avdual} \in \End^0(X \times X^\avdual)\Biggm| 
\vcenter{
\setbox0=\hbox{$\beta\in \Hom^{0,\sym}(X^\avdual,X)$}
\copy0
\hbox to\wd0{$\gamma\in \Hom^{0,\sym}(X,X^\avdual)$\hfill}
}\Biggr\}\, .
$$
Here $\Hom^{0,\sym}(X,X^\avdual)$ is the space of elements $\gamma \in  \Hom(X,X^\avdual) \otimes \mQ$ such that $\gamma=\gamma^\avdual$; similarly for $\Hom^{0,\sym}(X^\avdual,X)$.

In this section we construct an action of $\sp(X \times_S X^\avdual)$ on the motive $R(X/S)$ in the category $\Mot(S)$ of relative Chow motives with regard to ungraded correspondences; see Thm.~\ref{RefLefThm}. This connects with several ideas and results in the literature. On the one hand, it is related to the work of Looijenga and Lunts in~[\ref{LooijLunts}]. They introduced the idea that, given a projective variety~$Y$, rather than choosing one polarization and considering the resulting Lefschetz $\sl_2$-action on the cohomology, one may consider the Lie subalgebra $\gNS(Y)$ of $\gl\bigl(H^*(Y)\bigr)$ generated by the various $\sl_2$-triples obtained from all possible choices of a polarization. They obtained some general results about the structure of this ``N\'eron-Severi Lie algebra''~$\gNS(Y)$; in particular they proved (over~$\mC$) that it is semisimple. For $X$ a complex abelian variety they explicitly determined the structure of~$\gNS(X)$; it turns out that $\gNS(X)$ is an ideal of the Lie algebra that we call $\sp(X \times X^\avdual)$ and that in most cases the two are equal; see [\ref{LooijLunts}], Section~3. Our result provides a lifting of the tautological action of $\gNS(X)$ on the cohomology of~$X$ to an action on the Chow motive and it extends this to the full Lie algebra $\sp(X \times_S X^\avdual)$; a further generalization lies in the fact that we work in the setting of abelian schemes.

On the other hand, several people have obtained an action of an algebraic group~$U$, whose Lie algebra is, or is closely related to, $\sp(X \times_S X^\avdual)$ on, for instance, the derived category of~$X$. Results of this type have been obtained, independently, by Mukai~[\ref{MukaiSpin}], Polishchuk~[\ref{PolishThesis}], [\ref{Polish96}] and Orlov~[\ref{Orlov}]. (Note that various authors have used different names for the group~$U$ in question.) In [\ref{PolishWeilRep}] Polishchuk extended this to the setting of abelian schemes and obtained (loc.\ cit., Thm.~5.1) a kind of projective action on the relative motive. The Lie algebra action on the relative Chow motive that we construct has the advantage that it is a true action; in general it does not lift to an action of an algebraic group $\Sp(X\times X^\avdual)$.

\ssection
\ssectlabel{logGamma}
We start with some technical preparations. Let $p\colon A \to S$ be an abelian scheme over~$S$ of relative dimension~$g$, with zero section $e \colon S \to A$. Denote by $\CH^g_+(A/S)$ the kernel of the map $p_* \colon \CH^g(A) \to \CH^0(S)$, which is a nilpotent ideal of $\CH^g(A)$ for the $*$-product. (See [\ref{KimCorr}], Section~1.) For $y \in \bigl[e(S)\bigr] + \CH^g_+(A/S)$ we have a class $\log(y) \in \CH^g_+(A/S)$, defined by $\log(y) = \sum_{j\geq 1} (-1)^{j-1} y^j/j$; note that the sum is finite. See [\ref{KimCorr}], Def.~2.1 or also [\ref{KunChow}], (1.4.2).

If $a \in A(S)$ is a section, let $[\Gamma_a] = a_*[S]$; then it is clear that $[\Gamma_a] \in \bigl[e(S)\bigr] + \CH^g_+(A/S)$, so $\log[\Gamma_a]$ is defined. If 
$$
\def\normalbaselines{\baselineskip20pt\lineskip3pt\lineskiplimit3pt}
\matrix{A^\prime & \mapright{h} & A\cr
\mapdownl{p^\prime} && \mapdownr{p}\cr
S^\prime & \mapright{} & S\cr}
$$
is a cartesian square and $a^\prime \in A^\prime(S^\prime)$ is the section induced by~$a$, we have the relation $\log[\Gamma_{a^\prime}] = h^*\log[\Gamma_a]$.

The section~$a$ gives a map $i_a \colon A^\avdual \to A \times A^\avdual$, by $y \mapsto (a,y)$. If $\ell_A$ is the first Chern class of the Poincar\'e bundle on $A \times_S A^\avdual$, we have $\log[\Gamma_a] = (-1)^{g+1} \Four_{A^\avdual}(i_a^* \ell_A)$.

Consider an abelian scheme $X \to S$. We may view $X \times_S X$ as an abelian scheme over~$X$ via the first projection. The sections of $\pr_1 \colon X^2 \to X$ are in bijective correspondence with the endomorphisms of $X/S$. For an endomorphism $f \colon X \to X$ the associated class~$[\Gamma_f]$ is the usual class of the graph of~$f$. More generally, if $\alpha \in D = \End^0(X/S)$, write $\alpha = f/m$ for some $f \in \End(X/S)$ and $m \in \mZ\setminus\{0\}$; then we define 
$$
[\Gamma_\alpha] = \sum_{s \geq 0}\, {[\Gamma_f]_{(s)}\over m^s}\, ,
$$ 
where $[\Gamma_f]_{(s)}$ is the component of~$[\Gamma_f]$ in $\CH^g_{(s)}(X) = \CH^g\bigl(R^{2g-s}(X/S)\bigr)$. Again we have $[\Gamma_\alpha] \in \bigl[e(S)\bigr] + \CH^g_+(A/S)$; hence $\log[\Gamma_\alpha]$ is defined. The map $D \to \CH^g(X \times_S X)$ given by $\alpha \mapsto \log[\Gamma_\alpha]$ is $\mQ$-linear.

\ssection
\ssectlabel{logonProduct}
{\it Lemma. --- Let $S$ and~$T$ be connected schemes that are smooth of finite type over a Dedekind ring~$F$. Let $A \to S$ and $B\to T$ be abelian schemes, with zero sections $e_A$ and~$e_B$, respectively. Let $a \in A(S)$ and $b \in B(T)$ be sections, and consider the resulting section $(a\times b)$ of the abelian scheme $A \times_F B$ over $S \times_F T$. Then
$$
\log[\Gamma_{(a\times b)}] = \bigl(\log[\Gamma_a] \times [\Gamma_{e_B}]\bigr) + \bigl([\Gamma_{e_A}] \times \log[\Gamma_b]\bigr)
$$
as classes on $A \times_F B$.}
\medskip

\Proof
Let $g = \dim(A/S)$ and $h = \dim(B/T)$, and let $U = S\times_F T$. Write $A^\prime = U \times_S A$ and $B^\prime = U \times_T B$, and note that $A \times_F B \cong A^\prime \times_U B^\prime$ as $U$-schemes. We view $a$ (resp.~$b$) as sections of $A^\prime$ (resp.~$B^\prime$) over~$U$.  Let $\incl_{A^\prime} \colon A^\prime \to A^\prime \times_U B^\prime$ and $\incl_{B^\prime} \colon B^\prime \to A^\prime \times_U B^\prime$ be the maps given by $a \mapsto (a,0)$ and $b \mapsto (0,b)$. For the first Chern classes of the Poincar\'e bundles we have the relation
$$
\ell_{A^\prime\times_U B^\prime} = \pr^*_{A^\prime \times_U A^{\prime,\avdual}}(\ell_{A^\prime}) + \pr^*_{B^\prime \times_U B^{\prime,\avdual}}(\ell_{B^\prime})\, .
$$
This gives $i_{(a\times b)}^*(\ell_{A^\prime\times_U B^\prime}) = \pr_{A^{\prime,\avdual}}^* i_a^*(\ell_{A^\prime}) + \pr_{B^{\prime,\avdual}}^* i_b^*(\ell_{B^\prime})$. By [\ref{DPCRIFT}], (3.7.1), it follows that
$$
\Four_{A^{\prime,\avdual}\times_U B^{\prime,\avdual}}\bigl(i_{(a\times b)}^*(\ell_{A^\prime\times_U B^\prime})\bigr) =
(-1)^h \cdot \incl_{A^\prime,*} \Four_{A^{\prime,\avdual}}(i_a^*\ell_{A^\prime}) + (-1)^g \cdot \incl_{B^\prime,*} \Four_{B^{\prime,\avdual}}(i_b^*\ell_{B^\prime})\, .
$$
Hence,
$$
\eqalign{\log[\Gamma_{(a\times b)}] &= (-1)^{g+1} \cdot \incl_{A^\prime,*} \Four_{A^{\prime,\avdual}}(i_a^*\ell_{A^\prime}) + (-1)^{h+1} \cdot \incl_{B^\prime,*} \Four_{B^{\prime,\avdual}}(i_b^*\ell_{B^\prime})\cr
&= \incl_{A^\prime,*} \bigl(\log[\Gamma_a]\bigr) + \incl_{B^\prime,*}\bigl(\log[\Gamma_b]\bigr)\, ,\cr}
$$
which is just $\bigl(\log[\Gamma_a] \times [\Gamma_{e_B}]\bigr) + \bigl([\Gamma_{e_A}] \times \log[\Gamma_b]\bigr)$.
\QED

\ssection
\ssectlabel{halphaDef}
We now define the operators that will give us a Lie algebra action on the motive. First we consider the endomorphisms of~$X$. As the action of $D^{\opp,*}$ on $R(X/S)$ is algebraic, it induces an action of the Lie algebra~$\gl(X)^\opp$. We convert this into an action of~$\gl(X)$. Explicitly, for $\alpha \in \gl(X)$, let $h^\sharp_\alpha$ denote the endomorphism of~$R(X/S)$ that acts on $R^i(X/S) = \wedge^i R^1(X/S)$ as
$$
-\bigl([\transp\Gamma_\alpha] \wedge \id \wedge \cdots \wedge \id\bigr) - \bigl(\id \wedge [\transp\Gamma_\alpha] \wedge \id \wedge \cdots \wedge \id\bigr) - \cdots - \bigl(\id \wedge \cdots \wedge \id \wedge [\transp\Gamma_\alpha]\bigr)\, . 
$$
(So, $\alpha \in D=\gl(X)$ acts on $R^1(X/S)$ as $-[\transp\Gamma_\alpha]$; the minus signs come from the isomorphism $D^* \isomarrow D^{\opp,*}$ given by $\alpha \mapsto \alpha^{-1}$. Note that the cohomological grading on~$R(X/S)$ corresponds to the operator $h^\sharp_{-\id_X} = -h^\sharp_{\id_X}$.)

An element of~$D$ has a characteristic polynomial (of degree~$2g$) and in particular also a trace. Let $\trace\colon \gl(X) \to \mQ$ be the trace map. (So $\trace(\id)=2g$ and if $D$ is simple $\trace$ is a multiple of the reduced trace over~$\mQ$.) Now define operators $h_\alpha \in \End_{\Mot^0(S)}\bigl(R(X/S)\bigr)$, for $\alpha \in \gl(X)$, by
$$
h_\alpha = h^\sharp_\alpha + {\trace(\alpha)\over 2} \cdot \id_{R(X/S)}\, .
$$
If the context requires it we write $h^\sharp_{X,\alpha}$ and~$h_{X,\alpha}$.

For $\alpha \in D^*$ we have the relation $\Four_X \circ [\transp\Gamma_\alpha] = [\Gamma_{\alpha^\avdual}] \circ \Four_X = \deg(\alpha) \cdot [\transp\Gamma_{\alpha^{-1,\avdual}}] \circ \Four_X$. (Cf.\ the proof of Prop.~\ref{FDualEffect}.) Taking derivatives we obtain $-\Four_X \circ h^\sharp_{X,\alpha} = \bigl(\trace(\alpha) + h^\sharp_{X^\avdual,\alpha^\avdual}\bigr) \circ \Four_X$, which gives
$$
\Four_X \circ h_{X,\alpha} = -h_{X^\avdual,\alpha^\avdual} \circ \Four_X \eqlabel{Fh=-hF}
$$

\ssection
\ssectlabel{egamfbetDef}
Next we define Lefschetz and Lambda operators as in~[\ref{KunLef}]. For $\gamma \colon X \to X^\avdual$ we define a class $\ell(\gamma) \in \CH^1\bigl(R^2(X/S)\bigr)$ by
$$
\ell(\gamma) = {1\over 2} (\id_X,\gamma)^* \ell\, ,
$$
half the pull-back of the Poincar\'e class~$\ell$ under $(\id_X,\gamma) \colon X \to X \times_S X^\avdual$. The map $\gamma \mapsto \ell(\gamma)$ is linear. For arbitrary $\gamma \in \Hom^0(X,X^\avdual)$, write $\gamma = \gamma^\prime/m$ with $\gamma^\prime$ a true homomorphism and $m$ a nonzero integer; then set $\ell(\gamma) = \ell(\gamma^\prime)/m$. Note that we do not assume $\gamma$ to be symmetric. If $\gamma = \gamma^\avdual$ then $\ell(\gamma) \in \CH^1\bigl(R^2(X/S)\bigr)$ and the map $\Hom^{0,\sym}(X,X^\avdual) \to \CH^1\bigl(R^2(X/S)\bigr)$ thus obtained is bijective.

Define endomorphisms $L_\gamma \in \End_{\Mot(S)}\bigl(R(X/S)\bigr)$, for $\gamma \in \Hom^0(X,X^\avdual)$, by 
$$
L_\gamma = \Delta_*\bigl(\ell(\gamma)\bigr) \in \CH^{g+1}(X \times_S X)\, .
$$
On Chow groups, $L_\gamma$ is the intersection product with the class~$\ell(\gamma)$. It is an endomorphism of degree~$+2$, by which we mean that it is the sum of morphisms $L_\gamma \colon R^i(X/S) \to R^{i+2}(X/S)\bigl(1\bigr)$. If the context requires it, we write $\ell_X(\gamma)$ and $L_{X,\gamma}$ to indicate on which abelian scheme we work.

Dually, for $\beta \in \Hom^0(X^\avdual,X)$ we define $\lambda(\beta) \in \CH^{g-1}\bigl(R^{2g-2}(X/S)\bigr)$ by 
$$
\lambda(\beta) = (-1)^{g+1} \Four_{X^\avdual}\bigl(\ell_{X^\avdual}(\beta)\bigr)\, .
$$
Next define an endomorphism $\Lambda_\beta \in \End_{\Mot(S)}\bigl(R(X/S)\bigr)$ of degree~$-2$ by the commutativity of the diagrams 
$$
\def\normalbaselines{\baselineskip20pt\lineskip3pt\lineskiplimit3pt}
\matrix{R^{i+2}(X/S)\bigl(1\bigr) & \sizedmapright{-L_{X^\avdual,\beta}}{\Lambda_{X,\beta}} & R^i(X/S)\cr
\mapdownlr{\Four}{\wr} && \mapdownlr{\wr}{\Four}\cr
R^{2g-2-i}(X^\avdual/S)\bigl(g-i-1\bigr) & \mapright{-L_{X^\avdual,\beta}} & R^{2g-i}(X/S)\bigl(g-i\bigr)\cr} \eqlabel{FLamLFdiag}
$$
Using the basic properties of the Fourier transform, one readily verifies that on Chow groups, $\Lambda_\beta$ is the $*$-product with the class~$\lambda(\beta)$.

\ssection
\ssectlabel{efExa}
{\it Example.\/} --- Suppose we have a symmetric relatively ample bundle~$M$ on $X/S$ with rigidification along the zero section. If $\gamma\colon X \to X^\avdual$ is the associated polarization, our class~$\ell(\gamma)$ is just $c_1(M)$. In [\ref{KunLef}] this class is called~$d$. Our operator~$L_\gamma$ is K\"unnemann's Lefschetz operator~$L = L_d$. Further, our class $\lambda(\gamma^{-1})$ is K\"unnemann's ``curve class''~$c$ and our~$\Lambda_{\gamma^{-1}}$ is his Lambda-operator ${}^c\Lambda = {}^c\Lambda_c$. In particular, [\ref{KunLef}], Thm.~3.3, gives the commutation relation 
$$
[\Lambda_{\gamma^{-1}},L_\gamma] = \sum_{i=0}^{2g} (g-i)\pi_i = h_{\id_X}\, . \eqlabel{KunThm3.3}
$$
This is a special case of the relations we prove in Thm.~\ref{RefLefThm} below. 

The commutation relation~\eqref{KunThm3.3} in fact holds for all quasi-isogenies~$\gamma$ in $\Hom^{0,\sym}(X,X^\avdual)$. The proof in~[\ref{KunLef}] is entirely based on formal identities that hold in this generality; the positivity of the class~$\ell(\gamma)$ plays no role. Alternatively, we may remark that the polarization classes lie Zariski-dense in $\Hom^{0,\sym}(X,X^\avdual)$ and that \eqref{KunThm3.3} defines a Zariski-closed subset of the open set of quasi-isogenies in $\Hom^{0,\sym}(X,X^\avdual)$.
\medskip

The rest of this section is devoted to a proof of the following result.

\ssection
\ssectlabel{RefLefThm}
{\it Theorem. --- The map $\sp(X \times_S X^\avdual) \to \End_{\Mot(S)}\bigl(R(X/S)\bigr)$ defined by
$$
\pmatrix{\alpha & \beta \cr \gamma & -\alpha^\avdual} \mapsto \Lambda_\beta + h_\alpha + L_\gamma
$$
is a homomorphism of Lie algebras.}

\ssection
\ssectlabel{PfStep1}
To prove the theorem we have to establish the following commutation rules
$$
\matrix{
\hfill \hbox{a)\quad} & [h_{\alpha_1},h_{\alpha_2}]= h_{[\alpha_1,\alpha_2]}\hfill &\qquad& \hfill \hbox{d)\quad} & [h_\alpha,L_\gamma]= L_{-\alpha^\avdual\gamma -\gamma\alpha} \hfill\cr
\hfill \hbox{b)\quad} & [L_{\gamma_1},L_{\gamma_2}]= 0 \hfill &\qquad& \hfill \hbox{e)\quad} &  [h_\alpha,\Lambda_\beta]= \Lambda_{\alpha\beta+\beta\alpha^\avdual} \hfill\cr
\hfill \hbox{c)\quad} & [\Lambda_{\beta_1},\Lambda_{\beta_2}] = 0\hfill &\qquad& \hfill \hbox{f)\quad} & [\Lambda_\beta,L_\gamma] = h_{\beta\gamma}\hfill\cr
}
$$
(Here the $\beta_i$ and~$\gamma_i$ are symmetric.) The first three of these are clear from the definitions. Next we prove~d), which is equivalent to $[L_\gamma,h^\sharp_\alpha] = L_{\alpha^\avdual\gamma + \gamma\alpha}$. This will follow from the relation
$$
(\alpha^* \wedge \id + \id \wedge \alpha^*)\bigl(\ell(\gamma)\bigr) = \ell(\alpha^\avdual\gamma + \gamma\alpha)\, .\eqlabel{CommRel1}
$$
Indeed, by direct calculation we see that $[L_\gamma,h^\sharp_\alpha]$ acts on Chow groups as the intersection product with the left hand side of~\eqref{CommRel1}. Now consider the commutative diagram
$$
\def\normalbaselines{\baselineskip20pt\lineskip3pt\lineskiplimit3pt}
\matrix{
R^1(X/S) \otimes R^1(X^\avdual/S) & \maprightu{\id \otimes [\transp\Gamma_\gamma]} & R^1(X/S) \otimes R^1(X/S) & \maprightu{[\transp\Gamma_\Delta]} & R^2(X/S)\cr
\Big\Vert && \mapdown{\id \otimes [\transp\Gamma_\alpha]} && \mapdown{\id \wedge [\transp\Gamma_\alpha]} \cr
R^1(X/S) \otimes R^1(X^\avdual/S) & \mapright{\id \otimes [\transp\Gamma_{\gamma\alpha}]} & R^1(X/S) \otimes R^1(X/S) & \mapright{[\transp\Gamma_\Delta]} & R^2(X/S)\cr}
$$
As discussed in~\refn{MotFourDual} we have $\ell \in \CH^1\bigl(R^1(X/S) \otimes R^1(X^\avdual/S)\bigr)$. Taking images in $\CH^1\bigl(R^2(X/S)\bigr)$, the commutativity of the diagram gives the relation $(\id \wedge \alpha^*)\bigl(\ell(\gamma)\bigr) = \ell(\gamma\alpha)$. Similarly, we have a commutative diagram
$$
\def\normalbaselines{\baselineskip20pt\lineskip3pt\lineskiplimit3pt}
\matrix{
R^1(X/S) \otimes R^1(X^\avdual/S) & \maprightu{\id \otimes [\transp\Gamma_\gamma]} & R^1(X/S) \otimes R^1(X/S) & \maprightu{[\transp\Gamma_\Delta]} & R^2(X/S)\cr
\mapdown{[\transp\Gamma_\alpha] \otimes \id} && \mapdown{[\transp\Gamma_\alpha] \otimes \id} && \mapdown{[\transp\Gamma_\alpha] \wedge \id} \cr
R^1(X/S) \otimes R^1(X^\avdual/S) & \mapright{\id \otimes [\transp\Gamma_\gamma]} & R^1(X/S) \otimes R^1(X/S) & \mapright{[\transp\Gamma_\Delta]} & R^2(X/S)\cr}
$$
As $(\alpha \otimes \id_{X^\avdual})^* \ell = (\id_X \otimes \alpha^\avdual)^* \ell$ this gives the relation $(\alpha^* \wedge \id)\bigl(\ell(\gamma)\bigr) = \ell(\alpha^\avdual\gamma)$. Together, this proves~\eqref{CommRel1}, which gives us the commutation relation~d). 

By duality, using~\eqref{Fh=-hF} and our definition of the operators~$\Lambda_\beta$, we obtain~e).

\ssection
\ssectlabel{ActionHomolMot}
The hardest part of the proof of Thm.~\ref{RefLefThm} is to establish the commutation relation~f) in~\ref{PfStep1}. Our proof of this relation is based on a refinement of K\"unnemann's calculations in~[\ref{KunLef}]. The main new ingredient is an expression for the operators~$h_\alpha$ as algebraic cycles in $X \times_S X$. We give such an expression in Prop.~\ref{halphacycle}. As a preparation, we first introduce some  related operators~$k_\alpha$ that are best described working from the homological perspective. 

For $i\geq 0$, set $\epsilon_i = \pi_{2g-i}$ and define $R_i(X/S) = R^i(X/S)^\vee = (X,\epsilon_i,g)$. The addition $\Sigma \colon X \times_S X \to X$ define morphisms $[\Gamma_\Sigma] \circ (\epsilon_i \otimes \epsilon_j) \colon R_i(X/S) \otimes R_j(X/S) \to R_{i+j}(X/S)$. On Chow groups this induces the $*$-product. These morphisms give $R_*(X/S) = \oplus_{i=0}^{2g} R_i(X/S)$ the structure of a graded $\mQ$-algebra in the category $\Mot^0(S)$. By iteration of this product we obtain maps $R_1(X/S)^{\otimes i} \to R_i(X/S)$ which restrict to isomorphisms $\wedge^i R_1(X/S) \isomarrow R_i(X/S)$.

For $\alpha \in \gl(X)$ we define an endomorphism $k_\alpha \in \End_{\Mot^0(S)}\bigl(R_*(X/S)\bigr) = \CH^g(X \times_S X)$ by
$$
k_\alpha = [\Gamma_\id] *_{\pr_1} \log[\Gamma_\alpha]
$$
where $*_{\pr_1}$ is the $*$-product on $X \times_S X$ relative to the first projection.

\ssection
\ssectlabel{kalphaProp}
{\it Proposition. --- Let $\alpha \in \gl(X)$.

{\rm (\romno1)} For all $i \geq 0$ we have $\epsilon_i \circ k_\alpha = k_\alpha \circ \epsilon_i$.

{\rm (\romno2)} We have $k_\alpha \circ \epsilon_1 = [\Gamma_\alpha] \circ \epsilon_1$.

{\rm (\romno3)} The diagram
$$
\def\normalbaselines{\baselineskip20pt\lineskip3pt\lineskiplimit3pt}
\matrix{
R_*(X/S) \otimes R_*(X/S) & \mapright{[\Gamma_\Sigma]} & R_*(X/S)\cr
\mapdownl{k_\alpha \otimes \id + \id \otimes k_\alpha} && \mapdownr{k_\alpha}\cr
R_*(X/S) \otimes R_*(X/S) & \mapright{[\Gamma_\Sigma]} & R_*(X/S)\cr}
$$
is commutative.}
\medskip

\Proof
By definition,
$$
\log[\Gamma_\alpha] = \sum_{j\geq 1} {(-1)^{j-1}\over j}\, \bigl([\Gamma_\alpha]-[\Gamma_e]\bigr)^{* j} = \sum_{j\geq 1} \sum_{m=0}^j\, {(-1)^{m-1}\over j} {j\choose m} [\Gamma_{m\alpha}]\, ,
$$
where all $*$-products are taken relative the first projection. (Note that the sums over $j\geq 1$ are finite.) Part~(\romno1) follows because $[\Gamma_\id] *_{\pr_1} [\Gamma_{m\alpha}] = [\Gamma_{\id+m\alpha}]$ and $\epsilon_i$ commutes with~$[\Gamma_\beta]$ for all $\beta\in\gl(X)$.

(\romno2) We have $\epsilon_1 = \log[\Gamma_\id]$, so [\ref{KimCorr}], Lemma 2.2(\romno3) gives $[\Gamma_\alpha] \circ \epsilon_1 = \log[\Gamma_\alpha]$. On the other hand, for all $\alpha$, $\beta \in \gl(X)$ we have $\log[\Gamma_{\alpha +\beta}] = \log[\Gamma_\alpha] + \log[\Gamma_\beta]$; hence
$$
[\Gamma_{\alpha+\beta}] \circ \epsilon_1 = \bigl([\Gamma_\alpha] \circ \epsilon_1\bigr) + \bigl([\Gamma_\beta] \circ \epsilon_1\bigr)\, .
$$
(Taking duals, this proves that the map~$r^{(1)}$ of Proposition~\ref{risHomom} is indeed a homomorphism of algebras.) It follows that $\bigl([\Gamma_\id] *_{\pr_1} [\Gamma_\beta]\bigr) \circ \epsilon_1 = [\Gamma_{\id+\beta}] \circ \epsilon_1 = [\Gamma_\id] \circ \epsilon_1 + [\Gamma_\beta] \circ \epsilon_1$ for all $\beta \in \gl(X)$. Hence,
$$
k_\alpha \circ \epsilon_1 = \left(\sum_{j\geq 1} \sum_{m=0}^j\, {(-1)^{m-1}\over j} {j\choose m}\right)\, [\Gamma_\id] \circ \epsilon_1 + \log[\Gamma_\alpha] \circ \epsilon_1 \, ,
$$
which gives the desired conclusion because $\sum_{m=0}^j\, {(-1)^{m-1}\over j} {j\choose m} = 0$ for all $j \geq 1$.

For (\romno3) we start with the cartesian diagram
$$
\def\normalbaselines{\baselineskip20pt\lineskip3pt\lineskiplimit3pt}
\matrix{X^2 \times_S X & \mapright{\Sigma \times \id} & X \times_S X\cr
\mapdownl{\pr_1} && \mapdownr{\pr_1}\cr
X^2 & \sizedmapright{\Sigma \times \id}{\Sigma} & X\cr}
$$
where $X^2 = X \times_S X$. Using [\ref{DenMur}], Prop.~1.2.1, we have
$$
\eqalign{
k_\alpha \circ [\Gamma_\Sigma] = \bigl([\Gamma_\id] *_{\pr_1} \log[\Gamma_\alpha]\bigr) \circ [\Gamma_\Sigma] &= (\Sigma \times \id)^*\bigl([\Gamma_\id] *_{\pr_1} \log[\Gamma_\alpha]\bigr) \cr
&= (\Sigma \times \id)^* [\Gamma_\id] *_{\pr_1}  (\Sigma \times \id)^*\log[\Gamma_\alpha] \cr
&= [\Gamma_\Sigma] *_{\pr_1} \log[\Gamma_{\alpha \circ \Sigma}]\, .\cr}
$$
(In the last two expressions the $*$-product is taken relative to $\pr_1 \colon X^2 \times_S X \to X^2$.) By [\ref{KimCorr}], Lemma~2.2(\romno3),
$$
\log[\Gamma_{\alpha \circ \Sigma}] = \log[\Gamma_{\Sigma \circ (\alpha \times \alpha)}] = [\Gamma_\Sigma] \circ \log[\Gamma_{(\alpha \times \alpha)}]\, ;
$$
hence,
$$
\eqalign{
[\Gamma_\Sigma] *_{\pr_1} \log[\Gamma_{\alpha \circ \Sigma}] &= \bigl([\Gamma_\Sigma] \circ [\Gamma_{\id_{X^2}}]\bigr) *_{\pr_1} \bigl([\Gamma_\Sigma] \circ \log[\Gamma_{(\alpha \times \alpha)}]\bigr)\cr
&= [\Gamma_\Sigma] \circ \Bigl([\Gamma_{\id_{X^2}}] *_{\pr_1} \log[\Gamma_{(\alpha \times \alpha)}]\Bigr)\, .\cr}
$$
(Note that [\ref{KimCorr}], Lemma~2.2(\romno3), only applies to true endomorphisms~$\alpha$ of~$X$; it is, however, easy to verify that the identity we use is valid for arbitrary $\alpha \in D = \End^0(X/S)$.) 

For classes $\xi$ and~$\eta$ on $X \times_S X$, define their exterior product (as correspondences) $\xi \sqtensor \eta$ on $X^2 \times_S X^2$ by $\xi \sqtensor \eta = \pr_{13}^*(\xi) \cdot \pr_{24}^*(\eta)$. Lemma~\refn{logonProduct} gives $\log[\Gamma_{(\alpha \times \alpha)}] = \log[\Gamma_\alpha] \sqtensor [\Gamma_e] + [\Gamma_e] \sqtensor \log[\Gamma_\alpha]$, and since $[\Gamma_{\id_{X^2}}] = [\Gamma_{\id_X}] \sqtensor [\Gamma_{\id_X}]$ we conclude that
$$
[\Gamma_{\id_{X^2}}] *_{\pr_1} \log[\Gamma_{(\alpha \times \alpha)}] = k_\alpha \sqtensor [\Gamma_\id] + [\Gamma_\id] \sqtensor k_\alpha\, .
$$
Putting everything together we obtain the commutativity of the diagram in~(\romno3).
\QED

\ssection
\ssectlabel{kalphaCor}
{\it Corollary. --- The map $\alpha \mapsto k_\alpha$ gives $R_*(X/S)$ the structure of a graded $\gl(X)$-module in the category~$\Mot^0(S)$ such that $[\Gamma_\Sigma] \colon R_*(X/S) \otimes R_*(X/S) \to R_*(X/S)$ is a homomorphism of $\gl(X)$-representations. Via the isomorphism $\wedge^i R_1(X/S) \isomarrow R_i(X/S)$, the operator~$k_\alpha$ acts on $\wedge^i R_1(X/S)$ as
$$
([\Gamma_\alpha] \wedge \id \wedge \cdots \wedge \id) + (\id \wedge [\Gamma_\alpha] \wedge \id \wedge \cdots \wedge \id) + \cdots + (\id \wedge \cdots \wedge \id \wedge [\Gamma_\alpha])\, . 
$$
\par}
\vskip-\lastskip\medskip

\ssection
\ssectlabel{halphacycle}
{\it Proposition. --- For $\alpha \in \gl(X)$ the class $h_\alpha \in \CH^g(X \times_S X)$ is given by
$$
h_\alpha = k_\alpha - {\trace(\alpha)\over 2} \cdot [\Gamma_\id] = \bigl([\Gamma_\id] *_{\pr_1} \log[\Gamma_\alpha]\bigr)  - {\trace(\alpha)\over 2} \cdot [\Gamma_\id]\, .
$$\par}
\vskip-\lastskip\medskip

\Proof
We claim that $\Four_X \circ k_{X,\alpha} = -h^\sharp_{X^\avdual,\alpha^\avdual} \circ \Four_X$; combining this with~\eqref{Fh=-hF} then gives the assertion. The map $\alpha \mapsto \Four_X^{-1} \circ h^\sharp_{X^\avdual,-\alpha^\avdual} \circ \Four_X$ gives $R_*(X/S)$ the structure of a graded $\gl(X)$-module in the category~$\Mot^0(S)$ that is compatible with the $*$-product $[\Gamma_\Sigma]$. By Corollary~\ref{kalphaCor} it therefore suffices to prove that $k_\alpha$ and $\Four_X^{-1} \circ h^\sharp_{X^\avdual,-\alpha^\avdual} \circ \Four_X$ are equal on $R_1(X/S)$. This is just the standard relation $\Four_X \circ [\Gamma_{X,\alpha}] = [\transp\Gamma_{X^\avdual,\alpha^\avdual}] \circ \Four_X$.
\QED
\medskip

As a final preparation for the proof of Thm.~\ref{RefLefThm} we need to calculate the intersections of classes $\ell(\gamma)$ and~$\lambda(\beta)$.

\ssection
\ssectlabel{ellgamlambet}
{\it Lemma. --- For $\beta\in \Hom^{0,\sym}(X^\avdual,X)$ and $\gamma\in \Hom^{0,\sym}(X,X^\avdual)$ we have 
$$
\ell(\gamma) \cdot \lambda(\beta) = {\trace(\beta\gamma)\over 2}\cdot \bigl[e(S)\bigr] 
$$
in $\CH^g(X\times_S X)$.}
\medskip

\Proof
Because $\ell(\gamma) \in \CH^1\bigl(R^2(X/S)\bigr)$ and $\lambda(\beta) \in \CH^{g-1}\bigl(R^{2g-2}(X/S)\bigr)$, the class $\ell(\gamma) \cdot \lambda(\beta)$ lies in $\CH^g\bigl(R^{2g}(X/S)\bigr) = \mQ \cdot \bigl[e(S)\bigr]$. It therefore suffices to prove the assertion in the case that $S$ is the spectrum of a field, which from now on we assume. Moreover, we only have to prove that $\ell(\gamma) \cdot \lambda(\beta)$ has the correct degree.

By [\ref{DPCRIFT}], (3.7.1),
$$
\lambda(\beta) = {(-1)^{g+1}\over 2} \Four_{X^\avdual}\bigl((\beta,\id_{X^\avdual})^*\ell\bigr) = -{\textstyle {1\over 2}}\, (\beta,\id_{X^\avdual})^\avdual_* \bigl(\Four_{X \times_S X^\avdual}(\ell)\bigr)\, .
$$
By [\ref{BeauvFour}], \S 3, Lemme~1, we have $\Four(e^\ell) = e^{-\ell}$; hence $\Four_{X \times_S X^\avdual}(\ell) = (-1)^{2g-1} \switch_*\bigl(\ell^{[2g-1]}\bigr)$, where $\switch \colon X^\avdual \times_S X \isomarrow X \times_S X^\avdual$ is the map reversing the factors and where we recall that $\ell^{[n]} = \ell^n/n!$. It follows that $\lambda(\beta)$ is the push-forward of $\ell^{[2g-1]}/2$ under the morphism $[\id_X,\beta]\colon X \times_S X^\avdual \to X$ given by $(x,\xi) \mapsto x+\beta(\xi)$. This gives
$$
\ell(\gamma) \cdot \lambda(\beta) = {\textstyle {1\over 4}}\, [\id,\beta]_*\bigl(\ell^{[2g-1]}\bigr) \cdot (\id,\gamma)^*\bigl(\ell\bigr) = {\textstyle {1\over 4}}\, [\id,\beta]_*\Bigl(\ell^{[2g-1]} \cdot [\id,\beta]^* (\id,\gamma)^*\bigl(\ell\bigr)\Bigr) \, .
$$
Note that $(\id,\gamma) \circ [\id,\beta] \colon X\times_S X^\avdual \to X \times_S X^\avdual$ is the map given by the matrix $\left({\id\atop \gamma} {\beta\atop \gamma\beta}\right)$. 

In general, let $M$ be a non-degenerate line bundle on an abelian variety~$Y$ of dimension~$n$, with associated isogeny $\phi_M \colon Y \to Y^\avdual$. If $f \mapsto f^\prime$ is the corresponding Rosati involution on $\End^0(Y)$, we have
$$
\int_Y c_1(M)^{[n-1]} \cdot f^*\bigl(c_1(M)\bigr) = {\deg(\phi_M)^{1/2}\over 2} \cdot \trace(ff^\prime)\, .
$$
See [\ref{MAV}], \S 21, Thm.~1. We apply this to the Poincar\'e bundle~$P$ on $X \times_S X^\avdual$ and note that $\phi_P = \switch$. Because $\beta=\beta^\avdual$ and $\gamma=\gamma^\avdual$, we obtain
$$
\eqalign{
\int_X [\id,\beta]_*\Bigl( \ell^{[2g-1]} \cdot {\textstyle \left({\id\atop \gamma} {\beta\atop \gamma\beta}\right)}^*\bigl(\ell\bigr)\Bigr) &= \int_{X\times_S X^\avdual} \ell^{[2g-1]} \cdot {\textstyle \left({\id\atop \gamma} {\beta\atop \gamma\beta}\right)}^*\bigl(\ell\bigr)\cr
&= {\textstyle {1\over 2}} \trace\left({\textstyle \left({\id\atop \gamma} {\beta\atop \gamma\beta}\right)} {\textstyle \left({\beta\gamma\atop \gamma} {\beta\atop \id}\right)} \right)
= {\textstyle {1\over 2}} \trace {\textstyle \left({2\beta\gamma\atop 2\gamma\beta\gamma} \, {2\beta\atop 2\gamma\beta}\right)} = 2\, \trace(\beta\gamma)\, ,\cr} 
$$
which gives the assertion.
\QED

\ssection
\ssectlabel{PfLastStep}
We now complete the proof of Thm.~\ref{RefLefThm}. For this we have to establish relation~f) in~\ref{PfStep1}. As the set of pairs $(\beta,\gamma)$ for which this holds is Zariski-dense in $\Hom^{0,\sym}(X^\avdual,X) \times \Hom^{0,\sym}(X,X^\avdual)$, we may assume $\beta$ is a quasi-polarization. Write $X^2 = X\times_S X$, let $\Sigma \colon X^2 \to X$ be the addition map and let $\pr_1$, $\pr_2 \colon X^2 \to X$ be the projections. As before we view~$X^2$ as an abelian scheme over~$X$ via~$\pr_1$. Let $q_1$, $q_2\colon X^2 \times_X X^2 \to X^2$ be the projections and $\Sigma^{(2)} \colon X^2 \times_X X^2 \to X^2$ the addition map.

By [\ref{KunLef}], Lemmas 1.1(\romno1) and~3.1(\romno1), the commutator $[\Lambda_\beta,L_\gamma]$ is given by the class
$$
\pr_2^*\bigl(\lambda(\beta)\bigr) *_{\pr_1} \bigl([\Gamma_\id] \cdot \pr_2^*(\ell(\gamma))\bigr) - 
\pr_2^*\bigl(\ell(\gamma)\bigr) \cdot \bigl([\Gamma_\id] *_{\pr_1} \pr_2^*(\lambda(\beta)) \bigr)\, ,
$$
which by [\ref{KunLef}], Lemma~3.2 equals
$$
\eqalign{
&- [\Gamma_\id] *_{\pr_1} \Bigl(\pr_2^*(\lambda(\beta)) \cdot \pr_2^*(\ell(\gamma))\Bigr)\cr
&\qquad - \Sigma^{(2)}_*\biggl\{ \Bigl([\Gamma_\id] \times \pr_2^*(\lambda(\beta))\Bigr) \cdot \Bigl(\Sigma^{(2),*}\pr_2^*(\ell(\gamma)) -q_1^*\pr_2^*(\ell(\gamma)) - q_2^*\pr_2^*(\ell(\gamma))\Bigr) \biggr\}\, .\cr}
\eqlabel{fe-efExpr}
$$
The first term is easy to calculate: By Lemma~\ref{ellgamlambet}, $\pr_2^*(\lambda(\beta)) \cdot \pr_2^*(\ell(\gamma)) = \trace(\beta\gamma)/2 \cdot \pr_2^*\bigl[e(S)\bigr]$; so,
$$
- [\Gamma_\id] *_{\pr_1} \Bigl(\pr_2^*(\lambda(\beta)) \cdot \pr_2^*(\ell(\gamma))\Bigr) = -{\trace(\beta\gamma)\over 2}
\cdot [\Gamma_\id]\, .
$$

To calculate the second term, we identify $X^2 \times_X X^2 = X^3$. Then~$\Sigma^{(2)}$ is $\id_X \times \Sigma \colon X^3 \to X^2$; further, $q_1 = \pr_{12}$ and $q_2=\pr_{13}$. So we can rewrite $\Sigma^{(2),*}\pr_2^*(\ell(\gamma)) -q_1^*\pr_2^*(\ell(\gamma)) - q_2^*\pr_2^*(\ell(\gamma))$ as $\pr_{23}^* \Sigma^*(\ell(\gamma)) - \pr_{23}^* \pr_1^*(\ell(\gamma)) - \pr_{23}^*\pr_2^*(\ell(\gamma))$, which is equal to $\pr_{23}^* (\id\times \gamma)^*\bigl(\ell\bigr)$. 

Next note that $[\Gamma_\id] \times \pr_2^*\bigl(\lambda(\beta)\bigr) = \pr_{12}^*[\Gamma_\id] \cdot \pr_{23}^*\bigl(\pr_2^*(\lambda(\beta))\bigr)$. The second term of~\eqref{fe-efExpr} therefore equals
$$
\eqalign{&-(\id_X \times \Sigma)_*\Bigl\{\pr_{12}^*[\Gamma_\id] \cdot \pr_{23}^*\bigl(\pr_2^*(\lambda(\beta)) \cdot (\id\times \gamma)^*(\ell) \bigr)\Bigr\}\cr
&\qquad = -(\id_X \times \Sigma)_*\Bigl\{\pr_{12}^*[\Gamma_\id] \cdot \pr_{13}^*\bigl(\pr_2^*(\lambda(\beta)) \cdot (\id\times \gamma)^*(\ell) \bigr)\Bigr\}\, ,\cr}
$$
which by definition of~$*_{\pr_1}$ equals $-[\Gamma_\id] *_{\pr_1} \bigl(\pr_2^*(\lambda(\beta)) \cdot (\id\times \gamma)^*(\ell) \bigr)$. Our assumption that $\beta$ is a quasi-isogeny, together with the symmetry of $\beta$ and~$\gamma$, allows us to write
$$
\eqalign{
(\id \times \gamma)^* \ell &= (\id \times \beta^{-1})^* (\id \times \gamma\beta)^* \ell\cr
&= (\id \times \beta^{-1})^* ((\gamma\beta)^\avdual \times \id)^* \ell\cr 
&= (\id \times \beta^{-1})^* (\beta\gamma \times \id)^* \ell\cr
&= (\beta\gamma \times \id)^* (\id \times \beta^{-1})^* \ell\, .\cr}
$$
Further, $\pr_2 = \pr_2 \circ (\beta\gamma \times \id)$; so, $\pr_2^*\bigl(\lambda(\beta)\bigr) = (\beta\gamma \times \id)^* \pr_2^*\bigl(\lambda(\beta)\bigr)$. Hence,
$$
\eqalignno{
\pr_2^*\bigl(\lambda(\beta)\bigr) \cdot (\id\times \gamma)^*\bigl(\ell\bigr) &= (\beta\gamma \times \id)^*\Bigl(\pr_2^*(\lambda(\beta)) \cdot (\id\times \beta^{-1})^*(\ell) \Bigr)\qquad\qquad\cr
&= -(\beta\gamma \times \id)^* \log[\Gamma_\id]\, , & \hbox{by [\ref{KunLef}], Thm.~2.3}\cr
&= -\log\bigl((\beta\gamma \times \id)^*[\Gamma_\id]\bigr)\cr
&= -\log\bigl[\Gamma_{\beta\gamma}\bigr]\, . \cr}
$$
In total this gives
$$
[\Lambda_\beta,L_\gamma] = -{\trace(\beta\gamma)\over 2}
\cdot [\Gamma_\id] + [\Gamma_\id] *_{\pr_1} \log[\Gamma_{\beta\gamma}] 
$$
which by Prop.~\ref{halphacycle} equals $h_{\beta\gamma}$. This finishes the proof of Thm.~\ref{RefLefThm}.
\QED

\section{The canonical generalized Lefschetz decomposition}{LefschDecomp}
\medskip

\noindent
The goal of this section is to deduce from the $\sp(X \times X^\avdual)$-action on the motive a canonical Lefschetz decomposition of $R(X/S)$ that gives a further refinement of the decomposition of Thm.~\ref{MotDecRefined}. The main results are Thm.~\ref{RefLefCor} that gives the Lefschetz decomposition and Thm.~\ref{LefCompStr} that describes the structure of the Lefschetz components.

Throughout, $X/S$ is an abelian scheme as in~\ref{MotBasics}, $D = \End^0(X/S)$ is the endomorphism algebra, and $\gg = \sp(X \times_S X^\avdual)$.

\ssection
\ssectlabel{ggReps}
The Lie algebra $\gg$ has a natural grading $\gg = \gg_{-2} \oplus \gg_0 \oplus \gg_2$, with
$$
\gg_{-2} = \Hom^{0,\sym}(X^\avdual,X)\, ,\quad
\gg_0 = \gl(X)\, ,\quad
\gg_2 = \Hom^{0,\sym}(X,X^\avdual)\, .
$$
For the universal enveloping algebra $\UEA = \UEA(\gg)$ we have 
$$
\UEA = \Sym(\gg_2) \cdot \UEA(\gg_0) \cdot \Sym(\gg_{-2})\, . \eqlabel{UEADec}
$$
Define $\hbar = \bigl({-1\atop 0}\, {0\atop -1}\bigr) \in \gg$, which corresponds to $-\id_X \in \gg_0$. The pair $(\gg,\hbar)$ is a Lefschetz pair in the sense of~[\ref{LooijLunts}]. (Note, however, that we work over~$\mQ$. Some results of~[\ref{LooijLunts}] are valid only over an algebraically closed field, even though this is not always stated.)

For a Lie algebra~$\cL$ over~$\mQ$, let $\Rep^\lfss(\cL)$ denote the category of representations of~$\cL$ on $\mQ$-vector spaces that are direct sums of finite dimensional irreducible representations. (Locally finite semisimple representations.) Let $\psi\colon \gg \to \gl(V)$ be in $\Rep^\lfss(\gg)$. We have a grading $V = \oplus_i V^i$, defined by the requirement that $\hbar$ acts on~$V^i$ as multiplication by~$i$. Each summand~$V^i$ is a representation of $\gg_0 = \gl(X)$. We define the primitive subspace of~$V$ by
$$
V^\prim = \bigl\{v \in V \bigm| \beta(v) = 0 \quad \hbox{for all $\beta \in \gg_{-2}$} \bigr\}\, ,
$$
which is again a $\gg_0$-module. It follows from Wallach's results in~[\ref{Wallach}] that $V^\prim$ is an object of $\Rep^\lfss(\gg_0)$ and that the functor $\Rep^\lfss(\gg) \to \Rep^\lfss(\gg_0)$ given by $V \mapsto V^\prim$ is fully faithful. (In loc.\ cit.\ only finite dimensional representations are considered but the extension to the categories $\Rep^\lfss$ is immediate.) In particular, a finite dimensional $\gg$-representation~$V$ is irreducible if and only if its primitive part $V^\prim$ is irreducible as a $\gg_0$-module.

If the representation~$\psi$ is irreducible, there is a non-negative integer~$m$ such that $V^\prim = V^{-m}$ and $V^i \neq 0$ if and only if $i \in \{-m,-m+2,\ldots,m-2,m\}$; following [\ref{LooijLunts}] we call this integer~$m$ the depth of~$\psi$, notation $\depth(\psi)$.

\ssection
\ssectlabel{RefLefCor}
{\it Theorem. --- Let $\Irrep(\gg)$ be the set of isomorphism classes of finite dimensional irreducible representations of $\gg = \sp(X \times_S X^\avdual)$. 

{\rm (\romno1)} There is a unique decomposition
$$
R(X/S) = \bigoplus_{\psi \in \Irrep(\gg)}\, R_\psi(X/S) \eqlabel{MotLefDec}
$$
in the category $\Mot(S)$ that is stable under the action of~$\gg$ and has the property that for any object $M \in \Mot(S)$ the $\gg$-representation $\Hom_{\Mot(S)}\bigl(M,R_\psi(X/S)\bigr)$ is $\psi$-isotypic. There are finitely many $\psi \in \Irrep(\gg)$ for which $R_\psi(X/S)$ is non-zero. 

{\rm (\romno2)} The decomposition\/~{\rm \eqref{MotLefDec}} is stable under the action of~$D^{\opp,*}$ and there is a unique decomposition
$$
R(X/S) = \bigoplus_{\psi \in \Irrep(\gg)\atop \xi \in \chargp^\adm/\Gamma} \, R_\psi^{(\xi)}(X/S)
$$
in $\Mot^0(S)$ that is a common refinement of the decompositions\/~{\rm \eqref{RX/SDec}} and\/~{\rm \eqref{MotLefDec}}.

{\rm (\romno3)} If $\Phi \colon \Mot(S) \to \Vect_\mQ$ is a $\mQ$-linear functor, $\Phi\bigl(R_\psi(X/S)\bigr)$ is $\psi$-isotypic as a representation of~$\gg$.}
\medskip

If, in~(\romno2), the endomorphism algebra~$D$ is not simple, we proceed as explained in Remark~\ref{GeneralCase}.
\medskip

\Proof~Let $\tau$ denote the action of the algebraic group~$D^*$ on $R(X/S)$ via the operators $[\transp\Gamma_{\alpha^{-1}}] = \alpha^{-1,*}$. The action of $\gg_0 = \gl(X)$ on~$R(X/S)$ that we consider is ${\rm d}\tau \otimes {\trace\over 2}$. Hence 
$$
\Hom_{\Mot(S)}\bigl(M,R(X/S)\bigr) = \bigoplus_{m\in \mZ\atop \xi \in \chargp^\adm/\Gamma}\, \Hom_{\Mot^0(S)}\bigl(M(m),R^{(\xi)}(X/S)\bigr)
$$ 
is a sum of finite dimensional $\gg_0$-modules, for any $M \in \Mot(S)$. Because the elements of $\gg_{-2}$ and $\gg_2$ act as nilpotent endomorphisms with index of nilpotency at most~$g$, it follows that $\Hom_{\Mot(S)}\bigl(M,R(X/S)\bigr)$ is a sum of finite dimensional $\gg$-modules. The Lie algebra $\gg$ is reductive and its center is contained in~$\gg_0$; see [\ref{LooijLunts}], Sect.~3. By [\ref{JacobsLie}], Sect.\ \Romno3.11, Thm.~17, it follows that the $\gg$-module $\Hom_{\Mot(S)}\bigl(M,R(X/S)\bigr)$ is completely reducible. Further, it is clear that in the isotypic decomposition
$$
\Hom_{\Mot(S)}\bigl(M,R(X/S)\bigr) = \bigoplus_{\psi \in \Irrep(\gg)} H_\psi\,  \eqlabel{HomMRgRep}
$$
there are only finitely many~$\psi$ that occur, independently of~$M$. 

By construction, \eqref{HomMRgRep} is stable under the action of~$\gg_0$, and as $D^{\opp,*}$ is a connected algebraic group it follows that this decomposition is stable under the action of~$D^{\opp,*}$, too. Hence we have a further decomposition
$$
\Hom_{\Mot(S)}\bigl(M,R(X/S)\bigr) = \bigoplus_{\psi \in \Lef(X)\atop \xi \in \chargp^\adm/\Gamma} H_\psi^{(\xi)} \eqlabel{HomMRgD*Rep}
$$
that refines both~\eqref{HomMRgRep} and the decomposition obtained from~\eqref{RX/SDec}.

As \eqref{HomMRgRep} and \eqref{HomMRgD*Rep} are clearly functorial in~$M$ (because $\gg$ and~$D^{\opp,*}$ act through the factor~$R(X/S)$), the Yoneda Lemma gives that they correspond to decompositions of the motive~$R(X/S)$. (Use that $\Mot(S)$ is pseudo-abelian.) This proves (\romno1) and~(\romno2). Part~(\romno3) is a formal consequence; the proof is the same as in~\ref{Dopp*Realis}, replacing the group algebra~$\mQ[D^{\opp,*}]$ by the universal enveloping algebra~$\UEA(\gg)$.
\QED

\ssection
\ssectlabel{LefCompFourier}
{\it Remark.\/} --- We have an isomorphism of Lie algebras
$$
\iota\colon \sp(X \times_S X^\avdual) \isomarrow \sp(X^\avdual \times_S X)
\qquad\hbox{given by}\qquad
\pmatrix{\alpha & \beta\cr \gamma & -\alpha^\avdual} \mapsto \pmatrix{-\alpha^\avdual & -\gamma\cr -\beta & \alpha}\, .
$$
It is immediate from \eqref{Fh=-hF} and the commutativity of diagram~\eqref{FLamLFdiag} that the Fourier transform restricts to isomorphisms
$$
\Four\colon R_\psi(X/S) \isomarrow R_{\psi \circ \iota^{-1}}(X^\avdual/S)\, .
$$

\ssection
We refer to the motives $R_\psi(X/S)$ as the {\it Lefschetz components\/} of~$R(X/S)$. We should like to analyse these components further.

Let $\psi \colon \gg \to \gl(V)$ be an irreducible representation of depth $m = \depth(\psi)$. We have
$$
R_\psi(X/S) = \bigoplus_{j=0}^m\, R_\psi^{g-m+2j}(X/S) = \bigoplus_{j=0}^m \bigoplus_{\xi \in \chargp^\adm/\Gamma \atop \wt{\xi} = g-m+2j}\, R_\psi^{(\xi)}(X/S)\, ,
$$
and we define the primitive part of $R_\psi^\prim(X/S)$ to be the lowest weight summand $R_\psi^{g-m}(X/S)$. It follows from the remarks in~\ref{ggReps} that there is a unique $\xi \in \chargp^\adm/\Gamma$ with $\wt{\xi} = g-m$ such that $R_\psi^\prim(X/S) = R_\psi^{(\xi)}(X/S)$.  

In the classical case, when $D=\mQ$, we have $\gg = \sl_{2,\mQ}$. If $\psi \colon \sl_2 \to \gl(V)$ is isomorphic to the $k$th symmetric power of the standard representation, the Lefschetz component $R_\psi(X/S)$ is isomorphic, in $\Mot(S)$, to $R_\psi^\prim(X/S) \otimes V$. In general the structure of the Lefschetz components is a little more subtle; for instance, it is not true that $R_\psi(X/S)$ is always a direct sum of copies of motives of the form $R_\psi^\prim(X/S)\bigl(-j\bigr)$. The problem is that the primitive part~$R_\psi^\prim$ is still ``too big'': it decomposes as a direct sum of copies of a smaller motive~$P_\psi$ that we call the {\it core\/} of the Lefschetz component. This core comes equipped with a right action of the division algebra $B = \End_\gg(\psi)$. The full Lefschetz component can then again be described as $P_\psi \otimes_B V$.

The next result makes this description precise. We consider a finite dimensional irreducible representation $\psi \colon \gg \to \gl(V)$. Let $\psi_0 \colon \gg_0 \to \gl(V^\prim)$ be the associated representation of~$\gg_0$ on the space of primitive vectors, and define $B = \End_\gg(V) = \End_{\gg_0}(V^\prim)$. (For the second equality, see~\ref{ggReps}.) Let $m = \depth(\psi)$, so that $V = \oplus_{j=0}^m\, V^{-m+2j}$ and $V^\prim = V^{-m}$.

\ssection
\ssectlabel{LefCompStr}
{\it Theorem. ---  With notation and assumptions as just described, the contravariant functor 
$$
P_\psi \colon \Mot^0(S) \to \Mod_B
$$
defined by $P_\psi(M) = \Hom_{\gg_0}\bigl(V^\prim,\Hom_{\Mot^0(S)}(M,R_\psi^\prim(X/S))\bigr)$ is representable by a motive~$P_\psi$ with right $B$-action and we have a natural isomorphism
$$
\bigoplus_{j=0}^m\, P_\psi(-j) \otimes_B V^{-m+2j} \isomarrow R_\psi(X/S) \eqlabel{CoreIsom}
$$
in $\Mot^0(S)$.}
\medskip

By $P_\psi(-j) \otimes_B V^{-m+2j}$ we mean the motive representing the functor 
$$
M \mapsto \Hom_{\Mot^0(S)}\bigl(M,P_\psi(-j)\bigr) \otimes_B V^{-m+2j}\, .
$$
As $V^{-m+2j}$ is a free $B$-module of finite rank, this motive is (non-canonically) isomorphic to a sum of copies of $P_\psi(-j)$. Note that in $\Mot(S)$ the LHS of \eqref{CoreIsom} can be written a little more suggestively as $P_\psi \otimes_B V$.
\medskip

\Proof~Write $R_\psi = R_\psi(X/S)$. With $m = \depth(\psi)$ we have $R_\psi = \oplus_{j=0}^m\, R_\psi^{g-m+2j}$ and $R_\psi^\prim = R_\psi^{g-m}$. Then $\oplus_{j=0}^m\, R_\psi^{g-m+2j}(j)$ is a $\psi$-isotypic object with $\gg$-action in~$\Mot^0(S)$.

As discussed in \ref{ggReps}, if $W$ is a $\psi$-isotypic $\gg$-module, the natural map $\Hom_\gg(V,W) \to \Hom_{\gg_0}(V^\prim,W^\prim)$ is an isomorphism. This gives us an isomorphism of contravariant functors $\Mot^0(S) \to \Rep^\lfss(\gg)$,
$$
P_\psi \otimes_B V \isomarrow \bigoplus_{j=0}^m\, R_\psi^{g-m+2j}(j)
$$
that restricts to isomorphisms $P_\psi \otimes_B V^{-m+2j} \isomarrow R_\psi^{g-m+2j}(j)$. In particular, if $r$ is the rank of~$V^\prim$ as a $B$-module we see that $R_\psi^\prim$ is isomorphic, as a functor, to a sum of $r$~copies of~$P_\psi$. So $\End_B(V^\prim) \cong M_r(B^\opp)$ acts on~$R_\psi^\prim$ and $P_\psi$ is representable by the sub-object of~$R_\psi^\prim$ that is cut out by a rank~$1$ idempotent in $\End_B(V^\prim)$. \QED

\ssection
\ssectlabel{HighWeightsRem}
{\it Remark.\/} --- Assume $D = \End^0(X/S)$ is a simple algebra with center~$K$. If ${\frak h} \subset \gg_0 = \gl(X)$ is a Cartan subalgebra, ${\frak h}$ is also a Cartan subalgebra of $\gg = \sp(X \times_S X^\avdual)$. In combination with Thm.~\ref{CHDecThm} this allows us to describe more precisely which $\gg$-representations occur in the Lefschetz decomposition~\eqref{MotLefDec}. The statement of the result is made slightly complicated by the fact that the action of~$\gg_0$ is not simply the derivative of the action of~$D^{\opp,*}$ considered in Section~\ref{MotivicDec}. To begin with, there is a sign coming in, because we consider an action of $\gl(X)$ rather than an action of~$\gl(X)^\opp$. Further we have to take into account the change from the ``naive'' operators~$h_\alpha^\sharp$ to the operators~$h_\alpha$; see~\ref{halphaDef}. With this in mind, let us now give the details.

Define 
$$
\tilde\Lambda^+ = \bigl\{(\lambda_1,\ldots,\lambda_d) \in \mQ^d \bigm|\lambda_i - \lambda_{i+1} \in \mZ_{\geq 0}\quad \hbox{for all $i=1,\ldots,d-1$} \bigr\}\, .
$$
and let $\tilde\chargp^+ = \oplus_{\sigma\in \Sigma(K)}\, \tilde\Lambda^+$. The irreducible (finite dimensional) representations of $\gg_0 = \gl(X)$ over~$\mQ$ are indexed by $\tilde\chargp^+/\Gamma$.

For $\lambda = (\lambda_1,\ldots,\lambda_d)$ and $q \in \mQ$, define $$
\check\lambda = (-\lambda_d,\ldots,-\lambda_1)
\qquad\hbox{and}\qquad
\lambda[q] = (q+\lambda_1,\ldots,q+\lambda_d)\, .
$$
We extend this notation to elements of~$\tilde\chargp^+/\Gamma$ in the usual way: if $\xi\in \chargp^+/\Gamma$ is the $\Gamma$-orbit of some $\blambda \in \chargp^+$, define $\check\xi$ and $\xi[q]$ in $\chargp^+_\mQ/\Gamma$ as the $\Gamma$-orbits of the elements $\check\blambda$ and~$\blambda[q]$ obtained by applying the operations~$\check{}$ and~$[q]$ componentwise. For instance, the element~$\xi^\wdual$ of~\ref{xiwdualDef} is $\check\xi\bigl[2g/nd\bigr]$. 

The connection between the action of $D^{\opp,*}$ on $R(X/S)$ and the action of the Lie subalgebra $\gg_0 \subset \gg$ is given by the rule that for $\xi \in \chargp^\adm/\Gamma$ the Lie algebra $\gg_0$ acts on the motive $R^{(\xi)}(X/S)$ through the representation $\check\xi\bigl[g/nd\bigr]$. As usual, this means that for any other motive~$M$ the $\gg_0$-representation $\Hom_{\Mot(S)}\bigl(M,R^{(\xi)}(X/S)\bigr)$ is $\check\xi\bigl[g/nd\bigr]$-isotypic. Note that in general $g/nd$ is only a half-integer, which is why we cannot expect that the action of $\gg = \sp(X \times_S X^\avdual)$ lifts to an action of an algebraic group. (If $D=\mQ$ this problem does not occur, and indeed, in that case we have an action of~$\SL_2$, see~[\ref{BeauvSL2}].)

We obtain from this strong restrictions on the Lefschetz components of $R(X/S)$: For $\psi \in \Irrep(\gg)$, a necessary condition for $R_\psi(X/S)$ to be non-zero is that all irreducible factors of $\psi|_{\gg_0}$ correspond to elements in $\chargp^\Lef/\Gamma$, where we define
$$
\Lambda^\Lef = \Bigl\{(\mu_1,\ldots,\mu_d) \in \bigl({g\over nd} + \mZ\bigr)^d \Bigm| {g\over nd} \geq \mu_1 \geq \cdots \geq \mu_d \geq -{g\over nd}\Bigr\}\, ,
$$
and let $\chargp^\Lef = \oplus_{\sigma\in \Sigma(K)}\, \Lambda^\Lef$. As we shall see in the examples discussed in the next section, in practice this gives us an easy method to determine a finite subset of $\Irrep(\gg)$ containing all~$\psi$ with $R_\psi(X/S) \neq 0$.

\section{Examples}{Examples}
\medskip

\noindent
In this section we illustrate our main results with two concrete examples.

\ssection
\ssectlabel{RealMultExa}
Consider an abelian scheme~$X/S$ whose endomorphism algebra is a totally real field~$K$ of degree~$n$ over~$\mQ$. In our previous notation this means we have $D=K$ and $d=1$. To make the example as concrete as possible, we further assume that the Galois group $\Gamma = \Gal(\Ktilde/\mQ)$ is the full symmetric group~$\gS_n$; so, $[\Ktilde:\mQ]=n!$. Let $h=g/n$, which is an integer. 

To index representations, consider the sets
$$
\Lambda^+ = \mZ \quad\supset\quad \Lambda^\pol = \mZ_{\geq 0} \quad\supset\quad \Lambda^\adm = \{0,1,\ldots,2h\}\, ,
$$
and for $? \in \{+,\pol,\adm\}$ define $\chargp^? = \oplus_{\sigma \in \Sigma(K)}\, \Lambda^?$.

Let $\ga_K$ be $K$ viewed as an abelian Lie algebra over~$\mQ$. For $\eta \in \chargp^+$ let $\tilde{U}^{[\eta]}$ denote the representation of~$\ga_K$ over~$\mQ$ with underlying space~$\Ktilde$, on which an element $a \in \ga_K = K$ acts as multiplication by $\sum_{\sigma \in \Sigma(K)} \eta(\sigma) \cdot \sigma(a)$. The isomorphism class of this representation only depends on the $\Gamma$-orbit $[\eta] \in \chargp^+/\Gamma$, which justifies the notation. The representation~$\tilde{U}^{[\eta]}$ is a sum of $\# \Stab_\Gamma(\eta)$ copies of an irreducible representation~$U^{[\eta]}$. To make this explicit, let $F(\eta) \subset \Ktilde$ be the subfield that corresponds to $\Stab(\eta) \subset \Gamma$; then it is clear that $F(\eta) \subset \Ktilde$ is stable under the action of~$\ga_K$ and $U^{[\eta]}$ is the representation $\ga_K \to \gl\bigl(F(\eta)\bigr)$ thus obtained. It follows from this description that $\End_{\ga_K}(U^{[\eta]})$ contains~$F(\eta)$. On the other hand, the image of the map $K \to F(\eta)$ given by $a \mapsto \sum_{\sigma \in \Sigma(K)} \eta(\sigma) \cdot \sigma(a)$ generates~$F(\eta)$ as a field (by Galois theory); so in fact $\End_{\ga_K}(U^{[\eta]}) = F(\eta)$.

Because $\Gamma = \gS_n$, to give a class $[\eta] \in \chargp^+/\Gamma$ is equivalent to giving the $n$ values taken by~$\eta$. For any unordered $n$-tuple of integers $[b_1,\ldots,b_n]$ we therefore have a well-defined irreducible representation $U^{[b_1,\ldots,b_n]}$ of~$\ga_K$. For instance, the $1$-dimensional representation defined by trace map $\trace_{K/\mQ} \colon K \to \mQ$ is the representation $U^{[1,\ldots,1]}$.

Consider the decomposition $R(X/S) = \oplus_{\xi \in \chargp^\adm/\Gamma} R^{(\xi)}$ of Theorem~\ref{MotDecRefined}. Let $\xi \in \chargp^\adm/\Gamma$ and write it as $\xi = [a_1,\ldots,a_n]$. Then for any motive~$M$ the $\gg_0$-representation $\Hom_{\Mot^0(S)}(M,R^{(\xi)})$ is a sum of copies of $U^{[h-a_1,\ldots,h-a_n]} = \trace_{K/\mQ}^{\otimes h} \otimes U^{[-a_1,\ldots,-a_n]}$; see Remark~\ref{HighWeightsRem}. Thus, the irreducible representations of~$\gg_0$ that occur in~$R(X/S)$ are the representations $U^{[b_1,\ldots,b_n]}$ for all unordered $n$-tuples $[b_1,\ldots,b_n]$ with $b_i \in \{-h,-h+1,\ldots,h\}$ for all~$i$.

The Lie algebra $\gg = \sp(X\times_S X^\avdual)$ is $\sl_{2,K}$ viewed as a Lie algebra over~$\mQ$. We have a grading $\gg = \gg_{-2} \oplus \gg_0 \oplus \gg_2$ with $\gg_0 = \ga_K$, and $\gg_{\pm 2} \cong U^{[\pm 2,0,\ldots,0]}$ as representations of~$\gg_0$. The finite dimensional irreducible representations of~$\gg$ are  parameterized by $\chargp^\pol/\Gamma$. The irreducible representation~$V_\psi$ corresponding to a $\Gamma$-orbit~$\psi$ in~$\chargp^\pol$ is a $\mQ$-form of the representation of $\gg \otimes_\mQ \Qbar = \oplus_{\sigma \in \Sigma(K)}\, \sl_{2,\Qbar}$ given by
$$
\bigoplus_{\mu \in \psi}\, \sqtensor\limits_{\sigma\in \Sigma(K)}\, \Sym^{\mu(\sigma)}(\Standard)\, ,
$$ 
with $\Standard$ the standard $2$-dimensional representation of $\sl_{2,\Qbar}$. As before, because we assume $\Gamma=\gS_n$ we may give~$\psi$ by an unordered $n$-tuple of integers, in which case we use the notation $V_{[b_1,\ldots,b_n]}$.

When we restrict the representation $V_{[b_1,\ldots,b_n]}$ to $\gg_0 \subset \gg$, the irreducible constituents that occur are those of the form $U^{[-b_1+2j_1,\ldots,-b_n+2j_n]}$ with $j_i \in \{0,\ldots,b_i\}$ for all~$i$. As these must all occur in~$R(X/S)$ it follows that $b_i \leq h$ for all~$i$. Hence the Lefschetz decomposition can be written as
$$
R(X/S) = \bigoplus_{0 \leq b_1 \leq \cdots \leq b_n \leq h}\, R_{[b_1,\ldots,b_n]}(X/S)\, .
$$

Let $\psi = [b_1,\ldots,b_n]$, which has depth $m = b_1+\cdots+b_n$. The primitive part of the representation~$V_\psi$ is the $\gg_0$-representation~$U^{[-b_1,\ldots,-b_n]}$. As discussed, the endomorphism algebra of this representation is isomorphic to the field~$F(\psi)$ that is obtained as the field of invariants in~$\Ktilde$ of the stabilizer of a representative of~$\psi$ in~$\chargp^\pol$. By what was discussed in~\ref{ggReps}, the $\gg$-representation~$V_\psi$ has the same endomorphism algebra. As $V_\psi^\prim$ is $1$-dimensional over~$F(\psi)$, we find that the core of the Lefschetz component $R_\psi(X/S)$ is its primitive part $R_\psi^{[-b_1,\ldots,-b_n]}(X/S)$, which is a motive with $F(\psi)$-action, and $R_\psi(X/S) = R_\psi^{[-b_1,\ldots,-b_n]}(X/S) \otimes_{F(\psi)} V_\psi$.

\ssection
\ssectlabel{DquatExaBis}
As our next example we consider the case where $D = \End^0(X/S)$ is a quaternion algebra over~$\mQ$ with $D \otimes_\mQ \mR \cong M_2(\mR)$. The relative dimension~$g$ is necessarily even; write $g = 2h$.

The irreducible representations of $\gg_0 = \bigl(D,[~,~]\bigr)$ over~$\mQ$ are indexed by the pairs $\lambda = (\lambda_1,\lambda_2)$ in~$\mQ^2$ with $\lambda_1 - \lambda_2 \in \mZ_{\geq 0}$. With $d(\lambda)$ as in~\eqref{dlambdaDef}, the corresponding representation $U^{(\lambda_1,\lambda_2)}$ is a $\mQ$-form of $d(\lambda)$ copies of $\Sym^{\lambda_1-\lambda_2}(\Standard) \otimes (\lambda_2 \cdot \trace)$, where $\Standard$ is the standard representation of~$\gl_2$. The endomorphism algebra $\End_{\gg_0}(U^{\lambda_1,\lambda_2})$ is $\mQ$ if $\lambda_1-\lambda_2$ is even and is isomorphic to~$D$ for $\lambda_1-\lambda_2$ odd.

As discussed in Example~\ref{DquatExa}, the motivic decomposition of Thm.~\ref{MotDecRefined} in this case takes the form $R(X/S) = \oplus_{g \geq \lambda_1 \geq \lambda_2 \geq 0}\, R^{(\lambda_1,\lambda_2)}$ with integral $\lambda_1$ and~$\lambda_2$. The Lie algebra~$\gg_0$ acts on the summand $R^{(\lambda_1,\lambda_2)}$ through the representation $U^{(h-\lambda_2,h-\lambda_1)}$; see Remark~\ref{HighWeightsRem}. Thus, the irreducible $\gg_0$-representations that occur in $R(X/S)$ are the representations $U^{(\mu_1,\mu_2)}$ with $(\mu_1,\mu_2) \in \mZ^2$ and $h \geq \mu_1 \geq \mu_2 \geq -h$.

The Lie algebra $\gg = \sp(X \times_S X^\avdual)$ is a non-split $\mQ$-form of $\sp_4$. The degree zero part~$\gg_0$ is a Levi subalgebra, and $\gg_2 \cong U^{(2,0)}$ and $\gg_{-2} \cong U^{(0,-2)}$ as representations of~$\gg_0$. 

The representation theory of~$\sp_4$ is described, with many examples, in [\ref{FultHar}], Chap.~16. As in loc.\ cit., \S 16.2, we denote, for non-negative integers $a$ and~$b$, by~$\Gamma_{a,b}$ the irreducible representation of~$\sp_4$ with highest weight $a\varpi_1 + b\varpi_2$, where $\varpi_1$ and $\varpi_2$ are the (short and long, respectively) fundamental dominant weights. (In the notation of loc.\ cit., $\varpi_1 = L_1$ and $\varpi_2 = L_1+L_2$.) 

The description of the irreducible representations of~$\gg$ over~$\mQ$ is not much different. Again we have irreducible representations~$V_{a,b}$ indexed by pairs $(a,b) \in \mZ_{\geq 0}^2$. The only difference with the split case is that $V_{a,b}$ is not, in general, absolutely irreducible; instead, $V_{a,b} \otimes \Qbar$ is a sum of $1$ or~$2$ copies of~$\Gamma_{a,b}$, depending on the parity of~$a$. 

To determine which $\gg_0$-representations occur in~$V_{a,b}$ we need to calculate the branching rule for the restriction of representations of~$\sp_4$ to its Levi subalgebra~$\gl_2$. While this is easy in concrete examples (see below), it seems cumbersome to give the exact rule for which $\gg_0$-representations occur in the general case. However, by looking at the weights it is easy to see that the primitive part of~$V_{a,b}$ is isomorphic to $U^{(-b,-a-b)}$, so a necessary condition for the $\gg$-representation~$V_{a,b}$ to occur in $R(X/S)$ is that $a+b\leq h$. (This is probably also a sufficient condition.) Hence, the Lefschetz decomposition becomes
$$
R(X/S) = \bigoplus_{0\leq a \leq a+b \leq h}\, R_{a,b}(X/S)
$$
with $R_{a,b}(X/S)$ the $V_{a,b}$-isotypical summand.

If $a$ is even, the primitive part $R_{a,b}^\prim(X/S)$ is isomorphic to $a+1$ copies of the core $P_{a,b}$, and $R_{a,b}(X/S) \cong P_{a,b} \otimes_\mQ V_{a,b}$ in $\Mot(S)$. If $a$ is odd, $R_{a,b}^\prim(X/S)$ is isomorphic to $(a+1)/2$ copies of the core $P_{a,b}$, which has a right action of~$D$, and $R_{a,b}(X/S) \cong P_{a,b} \otimes_D V_{a,b}$. Note that the Lefschetz components with $a$ even (resp.\ odd) occur in even (resp.\ odd) degrees.

Let us now specialize to some examples in low dimension.
\medskip

\noindent
\boxit{$g=2$} The Lefschetz decomposition is $R(X/S) = R_{0,0} \oplus R_{1,0} \oplus R_{0,1}$. The component $R_{1,0} = R^1(X/S) \oplus R^3(X/S)$ is the odd part of the motive; its core equals its primitive part~$R^1(X/S)$. 

For the even parts we have
$$
R_{0,0} = R^{(1,1)}
\qquad\hbox{and}\qquad
R_{0,1} = R^0(X/S) \oplus R^{(2,0)}(X/S) \oplus R^4(X/S)\, ,
$$ 
Over a field, $R^2(X) = R^{(2,0)} \oplus R^{(1,1)}$ is the decomposition of $R^2(X)$ into an algebraic and a transcendental part as constructed by Kahn, Murre and Pedrini in~[\ref{KMP}]. (Caution: in Hodge theory the algebraic part of the $H^2$ is purely of type $(1,1)$. Here it is the summand $R^{(2,0)}$. The upper indices $(\lambda_1,\lambda_2)$ only tell how $D^{\opp,*}$ acts; they have no direct relation with the Hodge bidegrees.)
\medskip

\noindent
\boxit{$g=4$} The Lefschetz decomposition is $R(X/S) = R_{0,0} \oplus R_{1,0} \oplus R_{0,1} \oplus R_{2,0} \oplus R_{1,1} \oplus R_{0,2}$. For the summands~$R_{a,b}$ in odd degrees ($a$~odd), calculation of the branching rules gives  
$$
R_{1,0} = R_{1,0}^{(2,1)} \oplus R_{1,0}^{(3,2)} \quad\hbox{and}\quad
R_{1,1} = R_{1,1}^{(1,0)} \oplus R_{1,1}^{(3,0)} \oplus R_{1,1}^{(2,1)} \oplus R_{1,1}^{(4,1)} \oplus R_{1,1}^{(3,2)} \oplus R_{1,1}^{(4,3)}\, .
$$
In both cases, the core is the primitive part. For the Lefschetz components in even degrees, we have
$$
\matrix{
R_{0,0} = R_{0,0}^{(2,2)}\hfill  &  R_{2,0} = R_{2,0}^{(2,0)} \oplus R_{2,0}^{(3,1)} \oplus R_{2,0}^{(2,2)} \oplus R_{2,0}^{(4,2)}\hfill \cr
\noalign{\vskip 4pt}
R_{0,1} = R_{0,1}^{(1,1)} \oplus R_{0,1}^{(3,1)} \oplus R_{0,1}^{(3,3)} \quad\hfill &  R_{0,2} = R_{0,2}^{(0,0)} \oplus R_{0,2}^{(2,0)} \oplus R_{0,2}^{(4,0)} \oplus R_{0,2}^{(2,2)} \oplus R_{0,2}^{(4,2)} \oplus R_{0,2}^{(4,4)}\cr}
$$
For instance, we see that $R_{0,2}$ is the purely algebraic part of the motive, which is a sum of Tate motives. At the other extreme, $R_{0,0}$ is the ``transcendental part'' of the middle degree motive~$R^4(X/S)$, i.e., the part that does not come from smaller degrees.

Let us further look at the summand $R_{2,0}$. The representation~$U_{2,0}$ is the adjoint representation. (The weight diagram can be found on page~246 of~[\ref{FultHar}].) The depth is~$2$; further, $\dim(U_{2,0}^{-2}) = 3 = \dim(U_{2,0}^2)$ and $\dim(U_{2,0}^0) = 4$. So $R_{2,0}$ is an example of a Lefschetz component with the property that $R_{2,0}^4$ is not isomorphic to a sum of copies of $R_{2,0}^\prim(-1) = R_{2,0}^2(-1)$. The core $P_{2,0} = P_{2,0}(X/S)$ is a motive of rank~$5$ (over~$\mC$ its Hodge realization is an irreducible weight~$2$ Hodge structure with Hodge numbers $1-3-1$), and $R_{2,0} \cong P_{2,0}^{\oplus 3} \oplus P_{2,0}(-1)^{\oplus 4} \oplus P_{2,0}(-2)^{\oplus 3}$.

\section{Application: A question of C.~Voisin}{Applic1}
\medskip

\noindent
In this section we apply Thm.~\ref{RefLefThm} to answer a question of C.~Voisin. We prove the analogue for abelian varieties of a conjecture of Beauville~[\ref{BeauvSplitting}] about the so-called {\it weak splitting property\/}; see also [\ref{VoisCHHyperk}] and [\ref{VoisCHDec}].

To explain Voisin's question, consider an abelian variety~$X$ over a field. The descending filtration on $\CH(X)$ associated to Beauville's grading is expected to have the properties conjectured by Beilinson-Bloch and Murre. In particular, $\CH^j_{(s)}(X)$ should be zero for $s<0$ (Beauville's conjecture~($\hbox{F}_p$) of~[\ref{BeauvFour}], \S 5) and the cycle class map should be injective on $\CH_{(0)}(X) := \oplus_j\, \CH^j_{(0)}(X)$ (see [\ref{MurConj}], Conjecture~D). Almost nothing seems known about this in general. Voisin's question is the weaker question whether the cycle class map is injective on the subalgebra of $\CH(X)$ generated by the classes in $\CH^1_{(0)}$. Cor.~\ref{QVoisinCor} gives a positive answer to this; it is in fact a special case of a more general result, Thm.~\ref{QVoisinThm}.

\ssection
\ssectlabel{QVoisPrep}
With $S$ as in~\ref{MotBasics}, consider an abelian scheme $\pi \colon X \to S$. We have an isomorphism $\pi^* \colon \CH(S) \isomarrow \CH\bigl(R^0(X/S)\bigr) \subset \CH(X)$. We view $\CH(X)$ as a $\CH(S)$-algebra for the intersection product. Write $1_X = [X]$ for the identity element. Recall that the Lie algebra $\gg = \sp(X\times_S X^\avdual)$ has a natural grading $\gg = \gg_{-2} \oplus \gg_0 \oplus \gg_2$ with $\gg_2 = \Hom^{0,\sym}(X,X^\avdual)$. Further recall that the map $\gamma \mapsto \ell(\gamma)$ of~\ref{egamfbetDef} gives an isomorphism $\Hom^{0,\sym}(X,X^\avdual) \isomarrow \CH^1\bigl(R^2(X/S)\bigr)$ and that the Lefschetz operator~$L_\gamma$ acts on~$\CH(X)$ as the intersection product with the class~$\ell(\gamma)$.

\ssection
\ssectlabel{Rvarpi}
{\it Lemma. --- {\rm (\romno1)} There is a unique $\varpi \in \Irrep(\gg)$ such that $R^0_\varpi(X/S) \neq 0$. 

{\rm (\romno2)} The $\varpi$-isotypic subspace $\CH\bigl(R_\varpi(X/S)\bigr) \subset \CH(X)$ is the $\CH(S)$-subalgebra of $\CH(X)$ generated by the classes in $\CH^1\bigl(R^2(X/S)\bigr)$.}
\medskip

We may think of $R_\varpi(X/S)$ as the ``purely algebraic'' part of $R(X/S)$, generated (relative to~$S$) by the symmetric divisor classes. For instance, in the example of~\ref{DquatExaBis} with $g=2$ (resp.\ $g=4$) it is the Lefschetz component $R_{0,1}$ (resp.\ $R_{0,2}$).
\medskip

\Proof (\romno1)~For any motive~$M$ the $\gg_0$-representation $\Hom_{\Mot(S)}\bigl(M,R^0(X/S)\bigr)$ is $(\trace/2)$-isotypic. Hence if $\psi\colon \gg \to \gl(V)$ is an irreducible representation with  $R^0_\psi(X/S) \neq 0$, the corresponding irreducible representation~$V^\prim$ of~$\gg_0$ is the representation $(\trace/2)$. But as discussed in~\ref{ggReps}, the isomorphism class of~$\psi$ is determined by the $\gg_0$-module~$V^\prim$.

(\romno2) The $\varpi$-isotypic subspace $\CH\bigl(R_\varpi(X/S)\bigr) \subset \CH(X)$ is the $\UEA(\gg_2)$-subspace of $\CH(X)$ generated by its primitive vectors. The assertion follows from the remark that $R_\varpi^\prim(X/S) = R^0(X/S)$ together with the facts recalled in~\ref{QVoisPrep}.  \QED

\ssection
\ssectlabel{QVoisinThm}
{\it Theorem. --- Let $\cC$ be a $\mQ$-linear category and let $\Phi \colon \Mot(S) \to \cC$ be a $\mQ$-linear functor. Write $\Phi(X) = \Phi\bigl(R(X/S)\bigr)$ and consider the natural map
$$
u\colon \CH(X) = \Hom_{\Mot(S)}\bigl(\unitmot,R(X/S)\bigr) \to \Hom_\cC\bigl(\Phi(\unitmot),\Phi(X)\bigr)\, .
$$
If $u$ is injective on $\CH(S)$, it is injective on the $\CH(S)$-subalgebra of~$\CH(X)$ generated by the classes in $\CH^1\bigl(R^2(X/S)\bigr)$.\/}
\medskip

\Proof~By the lemma, the $\CH(S)$-subalgebra of~$\CH(X)$ generated by  $\CH^1\bigl(R^2(X/S)\bigr)$ equals $\CH\bigl(R_\varpi(X/S)\bigr)$. By functoriality, the map~$u$ is a map of $\gg$-modules. Hence, if there is a non-zero $y \in \Ker(u) \cap \CH\bigl(R_\varpi(X/S)\bigr)$ then the kernel of~$u$ contains a $\gg$-submodule isomorphic to~$\varpi$. But any such submodule contains non-zero element of $\CH(S)$.
\QED
\medskip

We now specialize to an abelian variety over a field. We use the notation $\CH^j_{(s)}(X) = \CH^j\bigl(R^{2j-s}(X)\bigr)$. In particular, $\CH^1\bigl(R^2(X/S)\bigr)$ is the space $\CH^1_{(0)}(X)$ of symmetric divisor classes.

\ssection
\ssectlabel{QVoisinCor}
{\it Corollary. --- Let $X$ be an abelian variety over a field~$F$. Let $y \in \CH(X)$ be an element that can be written as $y = P(\ell_1,\ldots,\ell_r)$ for some polynomial $P \in \mQ[t_1,\ldots,t_r]$ and classes $\ell_i \in \CH^1_{(0)}(X)$. If $y$ is numerically trivial, $y=0$.\/}
\medskip

\Proof~Apply the theorem to the natural functor $\Phi \colon \Mot(F) \to \Mot(F)_\num$, in which case $u$ is the quotient map $\CH(X) \twoheadrightarrow \CH(X)/{\equiv}$. Now use that $\CH(F) \cap \Ker(u) = (0)$.
\QED
\medskip

With similar arguments we can also prove a result about the next ``layer'' in Beauville's decomposition. We restrict our attention to the case where the endomorphism algebra~$D$ is simple.

\ssection
\ssectlabel{Rvarrho}
{\it Lemma. --- Let $X/S$ be an abelian scheme such that $D = \End^0(X/S)$ is a simple algebra. Then there is a unique $\varrho \in \Irrep(\gg)$ such that $R^1_\varrho(X/S) \neq 0$.}
\medskip

\Proof~It is immediate from the definitions that there is a unique class $\xi \in \chargp^\adm/\Gamma$ with $\wt{\xi}=1$. If $\psi \colon \gg \to \gl(V)$ is an irreducible representations with $R_\psi^1(X/S) \neq 0$ then $\psi$ has depth $g-1$ and the primitive part $V^\prim = V^{1-g}$ is the $\gg_0$-module corresponding to $\check\xi\bigl[g/nd\bigr]$; see~\ref{HighWeightsRem}. As $\psi$ is determined by the associated $\gg_0$-module~$V^\prim$ (see~\ref{ggReps}), the lemma follows. \QED

\ssection
\ssectlabel{AJInjectThm}
{\it Theorem. --- Let $X$ be an abelian variety over a field~$F$ such that $\End^0(X)$ is a simple algebra. Let $\cD \subset \CH(X)$ be the $\mQ$-subalgebra generated by the classes in $\CH^1_{(0)}(X)$. 

{\rm (\romno1)} Suppose $F=\mC$. Let $y \in \CH^j(X)$ be a class that lies in the $\cD$-submodule of~$\CH(X)$ generated by the elements of $\CH^1_{(1)}(X) = \CH^1\bigl(R^1(X)\bigr)$. If the Abel-Jacobi class 
$$
\AJ(y) \in H^{2j-1}(X,\mC)/\bigl(\Fil^j+H^{2j-1}(X,\mQ)\bigr)
$$ 
of~$y$ is trivial then $y=0$.

{\rm (\romno2)} Suppose $F$ is finitely generated over its prime field, and let $\ell$ be a prime number different from $\char(F)$. Choose a separable closure~$F^\sep$ of~$F$ and let $\Gamma_F = \Gal(F^\sep/F)$ and $\bar{X} = X_{F^\sep}$. Let $y \in \CH^j(X) \otimes \Ql$ be a class that lies in the $(\cD \otimes \Ql)$-submodule of $\CH(X) \otimes \Ql$ generated by the elements of $\CH^1_{(1)}(X)$. If the $\ell$-adic Abel-Jacobi class 
$$
\AJ_\ell(y) \in H^1_\cont\bigl(\Gamma_F,H^{2j-1}(\bar{X},\Ql)\bigr)
$$ 
of~$y$ is trivial then $y=0$.}
\medskip

\Proof~We only give the proof of~(\romno2); the argument for~(\romno1) is essentially the same (but easier because we do not need to pass to $\Ql$-coefficients). We assume $y\neq 0$ and derive a contradiction.

Let $\varrho \colon \gg \to \gl(V)$ be the irreducible representation of Lemma~\ref{Rvarrho}. The depth of this representation is $g-1$. The representation $\varrho_{\Ql} \colon (\gg \otimes \Ql) \to \gl(V \otimes \Ql)$ is again completely reducible; let $\theta_1,\ldots,\theta_r$ be its irreducible constituents, with $\theta_i \colon (\gg \otimes \Ql) \to \gl(W_i)$. Each~$W_i$ has a grading given by the action of the element~$\hbar$ (as in~\ref{ggReps}) and $W_i^{1-g} \neq 0$ for all~$i$.

The $\Ql$-vector space $\CH_{(1),\ell} := \oplus_{q=1}^g\, \CH^q_{(1)}(X) \otimes \Ql$ is a representation of $\gg \otimes \Ql$. By assumption, $y$ lies in the subrepresentation generated by $\CH^1_{(1)}(X) = \CH^1\bigl(R^1(X)\bigr)$, and as $R^1(X/S) = R^1_\varrho(X/S)$, it follows that the $(\gg \otimes \Ql)$-subrepresentation of $\CH_{(1),\ell}$ generated by~$y$ is isomorphic to $\oplus_{i \in I}\, \theta_i^{\oplus m_i}$ for some subset $I \subset \{1,\ldots,r\}$ and multiplicities $m_i$. In particular, this submodule contains a non-zero element~$y_1$ that is homogeneous of degree $1-g$ (meaning that $\hbar(y_1) = (1-g) \cdot y_1$); for this element we have $y_1 \in \CH^1_{(1)}(X) \otimes \Ql$.

The kernel of the $\ell$-adic Abel-Jacobi map
$$
\AJ_\ell \colon \CH_{(1),\ell} \to \bigoplus_{q=1}^g\, H^1_\cont\bigl(\Gamma_F,H^{2q-1}(\bar{X},\Ql)\bigr)
$$ 
is a $(\gg \otimes \Ql)$-subrepresentation of $\CH_{(1),\ell}$. Hence if $\AJ_\ell(y) = 0$ then also $y_1 \in \Ker(\AJ_\ell)$. But $\AJ_\ell$ is injective on $\CH_{(1)}^1(X) \otimes \Ql$, so we arrive at the desired contradiction. \QED

\vskip2.0\bigskipamount plus 2pt minus 1pt%
\writetocentry{nnsect}{References}
\goodbreak\centerline{{\bf References}}%
\nobreak\vskip.75\bigskipamount plus 2pt minus 1pt%

{\eightpoint
\bibitem{Ancona}
G.~Ancona, {\it D\'ecomposition du motif d'un sch\'ema ab\'elien universel.\/} Thesis, Paris~\Romno{13}, 2012.

\bibitem{BeauvFour}
A.~Beauville, {\it Quelques remarques sur la transformation de Fourier dans l'anneau de Chow d'une vari\'et\'e ab\'elienne.\/} Algebraic geometry (Tokyo/Kyoto, 1982), 238--260. Lecture Notes in Math.\ 1016, Springer, Berlin, 1983. 

\bibitem{BeauvChow}
A.~Beauville, {\it Sur l'anneau de Chow d'une vari\'et\'e ab\'elienne.\/}
Math.\ Ann.\ 273 (1986), 647--651. 

\bibitem{BeauvSplitting}
A.~Beauvile, {\it On the splitting of the Bloch-Beilinson filtration.\/} Algebraic cycles and motives, Vol.~2, 38--53.
London Math. Soc. Lecture Note Ser.\ 344, Cambridge Univ.\ Press, Cambridge, 2007. 

\bibitem{BeauvSL2}
A.~Beauville, {\it The action of $\SL_2$ on abelian varieties.\/}
J.\ Ramanujan Math.\ Soc.\ 25 (2010), 253--263.

\bibitem{DelCatTens}
P.~Deligne, {\it Cat\'egories tensorielles.\/} Mosc.\ Math.\ J.\ 2 (2002), 227--248. 

\bibitem{DemGab}
M.~Demazure, P.~Gabriel, {\it Groupes alg\'ebriques. Tome I.\/} Masson \& Cie, Paris; North-Holland, Amsterdam, 1970.

\bibitem{DenMur}
C.~Deninger, J.~Murre, {\it Motivic decomposition of abelian schemes and the Fourier transform.\/} J.\ reine angew.\ Math.\ 422 (1991), 201--219. 

\bibitem{FultHar}
W.~Fulton, J.~Harris, {\it Representation theory. A first course.\/} Graduate Texts in Mathematics 129. Springer-Verlag, New York, 1991.

\bibitem{JacobsLie}
N.~Jacobson, {\it Lie algebras.\/}
Republication of the 1962 original. Dover Publications, Inc., New York, 1979.

\bibitem{KMP}
B.~Kahn, J.~Murre, C.~Pedrini, {\it On the transcendental part of the motive of a surface.\/} Algebraic cycles and motives, Vol.~2, 143--202.
London Math. Soc. Lecture Note Ser.\ 344, Cambridge Univ.\ Press, Cambridge, 2007.

\bibitem{KimCorr}
S.~Kimura, {\it Correspondences to abelian varieties.~I.\/}
Duke Math.\ J.\ 73 (1994), 583--591. 

\bibitem{Kings}
G.~Kings, {\it Higher regulators, Hilbert modular surfaces, and special values of L-functions.\/} Duke Math. J. 92 (1998), 61--127.

\bibitem{KunLef}
K.~K\"unnemann, {\it A Lefschetz decomposition for Chow motives of abelian schemes.\/}
Invent. Math. 113 (1993), 85--102. 

\bibitem{KunChow}
K.~K\"unnemann, {\it On the Chow motive of an abelian scheme.\/} Motives (Seattle, WA, 1991), Part~1, 189--205. Proc.\ Sympos.\ Pure Math.\ 55, Amer.\ Math.\ Soc., Providence, RI, 1994. 

\bibitem{KunArChow}
K.~K\"unnemann, {\it Arakelov Chow groups of abelian schemes, arithmetic Fourier transform, and analogues of the standard conjectures of Lefschetz type.\/} Math.\ Ann.\ 300 (1994), 365--392.

\bibitem{LooijLunts}
E.~Looijenga, V.~Lunts, {\it A Lie algebra attached to a projective variety.\/}
Invent.\ Math.\ 129 (1997), 361--412. 

\bibitem{DPCRIFT}
B.~Moonen, A.~Polishchuk, {\it Divided powers in Chow rings and integral Fourier transforms.\/} Adv.\ Math.\ 224 (2010), 2216--2236.


\bibitem{MukaiSpin}
S.~Mukai, {\it Abelian variety and spin representation.\/} Preprint, Warwick, 1998, 13pp. Available for download at http://www.kurims.kyoto-u.ac.jp/\kern-1pt\lower 3pt\hbox{\~{\ }}\kern-3pt mukai/

\bibitem{MAV}
D.~Mumford, {\it Abelian varieties.\/}
Tata Institute of Fundamental Research Studies in Math.\ 5. Oxford University Press, Oxford, 1970.

\bibitem{MurConj}
J.~Murre, {\it On a conjectural filtration on the Chow groups of an algebraic variety. I. The general conjectures and some examples.\/} 
Indag.\ Math.\ (N.S.) 4 (1993), 177--188. 

\bibitem{Orlov}
D.~Orlov, {\it Derived categories of coherent sheaves on abelian varieties and equivalences between them.\/} Izv.\ Ross.\ Akad.\ Nauk Ser.\ Mat.\ 66 (2002), 131--158; English translation in Izv.\ Math.\ 66 (2002), 569--594.

\bibitem{PolishThesis}
A.~Polishchuk, {\it Biextensions, Weil representation on derived categories and theta-functions.\/} Harvard thesis, 1996.

\bibitem{Polish96}
A.~Polishchuk, {\it Symplectic biextensions and a generalization of the Fourier-Mukai transform.\/} Math.\ Res.\ Lett.\ 3 (1996), 813--828.

\bibitem{PolishWeilRep}
A.~Polishchuk, {\it Analogue of Weil representation for abelian schemes.\/}
J. reine angew.\ Math.\ 543 (2002), 1--37. 

\bibitem{Tits}
J.~Tits, {\it Repr\'esentations lin\'eaires irr\'eductibles d'un groupe r\'eductif sur un corps quelconque.\/} J.\ reine angew.\ Math.\ 247 (1971), 196--220.

\bibitem{VoisCHHyperk}
C.~Voisin, {\it On the Chow ring of certain algebraic hyper-K\"ahler manifolds.\/} Pure Appl. Math. Q. 4 (2008), 613--649.

\bibitem{VoisCHDec}
C.~Voisin, {\it Chow rings and decomposition theorems for families of K3 surfaces and Calabi-Yau hypersurfaces.\/} Geom.\  Topol.\ 16 (2012), 433--473.

\bibitem{Wallach}
N.~Wallach, {\it Induced representations of Lie algebras and a theorem of Borel-Weil.\/} Trans.\ Amer.\ Math.\ Soc.\ 136 (1969), 181--187. 

\bibitem{WedhOrdin}
T.~Wedhorn, {\it Ordinariness in good reductions of Shimura varieties of PEL-type.\/} Ann.\ scient.\ \'Ec.\ Norm.\ Sup.\ (4) 32 (1999), 575--618.
\vskip 1cm

\noindent
University of Amsterdam, KdV Institute for Mathematics, PO Box 94248, 1090 GE Amsterdam, The Netherlands. 
\smallskip

\noindent
bmoonen@uva.nl
\par}

\bye